\documentclass[11pt]{article}

\usepackage[letterpaper,margin=1in]{geometry}
\usepackage[utf8]{inputenc}
\usepackage[T1]{fontenc}
\usepackage{newtxtext}
\usepackage{microtype}
\usepackage{graphicx}
\usepackage{subcaption}
\usepackage{booktabs}
\usepackage{multirow}
\usepackage{amsmath,amssymb,mathtools}
\usepackage{amsthm}
\usepackage{newtxmath}
\usepackage{algorithm}
\usepackage{algorithmic}

\usepackage{float}
\usepackage{placeins}
\usepackage{tikz}
\usetikzlibrary{arrows.meta,positioning,calc,decorations.pathmorphing}
\usepackage{enumitem}
\usepackage[acronym,section=section]{glossaries}
\usepackage[numbers,sort&compress]{natbib}
\usepackage[colorlinks=true,linkcolor=blue,citecolor=blue,urlcolor=blue]{hyperref}
\usepackage[capitalize,noabbrev]{cleveref}
\usepackage{url}
\usepackage[disable,textsize=tiny]{todonotes}

\newcommand{\nodelab}{}

\let\oldincludegraphics\includegraphics
\renewcommand{\includegraphics}[2][]{%
  \IfFileExists{#2}{\oldincludegraphics[#1]{#2}}{%
    \fbox{\parbox[c][0.22\textheight][c]{0.88\linewidth}{\centering Missing figure:\\ \texttt{\detokenize{#2}}}}%
  }%
}

\theoremstyle{plain}
\newtheorem{theorem}{Theorem}[section]

\theoremstyle{definition}

\theoremstyle{remark}

\newacronym{kan}{KAN}{Kolmogorov-Arnold network}
\newacronym{pinn}{PINN}{physics-informed neural network}
\newacronym{pde}{PDE}{partial differential equation}
\newacronym{sp}{SP}{Stefan problem}
\newacronym{eop}{EOP}{elastic obstacle problem}
\newacronym{plop}{PLOP}{\(p\)-Laplacian obstacle problem}
\newacronym{bvp}{BVP}{boundary value problem}
\newacronym{mlp}{MLP}{multi-layer perceptron}
\newacronym{resnet}{ResNet}{residual neural network}

\makeglossaries

\begin{document}

\title{Kolmogorov-Arnold Networks for Free-Boundary Partial Differential Equations}
\author{
Tan Phuong Dong Le%
\thanks{This work was conducted while the author was a Ph.D. candidate in the Department of Applied Mathematics, University of Waterloo, Canada. Email: \href{mailto:pdle@uwaterloo.ca}{pdle@uwaterloo.ca}. Corresponding author.}
\\[0.4em]
\small Department of Statistics and Actuarial Science, University of Waterloo, Canada
}

\date{}

\maketitle

\begin{abstract}
We study free-boundary problems within a physics-informed framework using Kolmogorov-Arnold network approximations. The proposed approach incorporates obstacle constraints, \gls{pde} inequalities, complementarity conditions, and boundary conditions through residual-based loss functions. We consider a linear elliptic obstacle problem, a nonlinear \(p\)-Laplacian obstacle problem, and a time-dependent one-phase Stefan problem. The proposed \gls{kan} solver is compared with \gls{pinn} and residual-network baselines. Numerical experiments show that \glspl{kan} achieve low relative \(L^2\) and \(L^\infty\) errors while accurately resolving contact regions and moving interfaces. The results indicate that \gls{kan} representations provide an effective alternative for solving free-boundary \glspl{pde}.
\end{abstract}


\vspace{1em}

\section{Introduction}
\label{sec:introduction}
Free-boundary \glspl{pde} arise when the unknown solution is coupled with a geometric feature such as a contact region, phase interface, or moving boundary \cite{ros2014extremal,ros2016nonlocal,ros2018obstacle}. These problems are more difficult to approximate than classical boundary-value problems because a numerical method must determine both the solution and the location of the free boundary \cite{ros2016nonlocal}. Motivated by recent advances in \gls{pinn} methods and neural representations \cite{raissi2019physics,sirignano2018dgm,wang2021deep}, we investigate \gls{kan} models for free-boundary \glspl{pde}. We propose a residual-based \gls{kan} framework that incorporates the governing \gls{pde}, inequality constraints, complementarity conditions, and initial or boundary conditions into the training objective. We evaluate the framework on elliptic obstacle problems and a time-dependent Stefan problem.

\subsection{Background and Motivation}
A key difficulty in \glspl{eop} is identifying the contact set, where the solution coincides with the obstacle function \cite{ros2018obstacle}. In \glspl{sp}, the solution evolves together with a moving phase boundary \cite{friedman1975one,hanzawa1981classical,damlamian1977some,gupta2017classical}. Therefore, a numerical method must approximate both the solution and the free-boundary geometry while satisfying the governing constraints and interface conditions \cite{ros2018obstacle,tran20151,wang2021deep}.
\subsection{Related Work}
Classical methods for free-\glspl{bvp} include finite-difference and finite-element discretizations, as well as penalty and regularization techniques \cite{tran20151,burman2017galerkin,johnson1992adaptive}. These methods are well established but may require mesh refinement and explicit interface localization \cite{burman2017galerkin}.

\Glspl{pinn} commonly use \glspl{mlp} as solution representations \cite{raissi2019physics,sirignano2018dgm}. For free-boundary problems, the loss must also encode inequality constraints, complementarity conditions, and the geometry of the contact or moving-interface region \cite{wang2021deep}.

The \gls{kan} architecture was recently introduced as an alternative to standard \glspl{mlp} \cite{liu2024kan}. A \gls{kan} replaces fixed node-wise activation functions with learnable univariate functions on network edges, typically parameterized by splines. This structure provides local flexibility that may be useful for accurate interface resolution. In contrast, \glspl{mlp} use fixed nonlinear activation functions at the nodes and train linear weights between layers. In this work, we combine residual-based physics-informed learning with \gls{kan} representations for free-boundary \glspl{pde}.
\subsection{Contributions}
The main contributions of this work are as follows:
\begin{itemize}
\item We develop a \gls{kan}-based physics-informed framework for free-boundary \glspl{pde}.

\item We construct residual-based losses that enforce the governing \gls{pde}, boundary conditions, obstacle constraints, complementarity conditions, and moving-interface conditions.

\item We evaluate the framework on the \gls{eop}, \gls{plop}, and \gls{sp}, and compare its accuracy and parameter efficiency with \gls{pinn} and \gls{resnet} baselines.
\end{itemize}
The remainder of the paper is organized as follows. Section~\ref{sec:problem-statement} introduces the \gls{kan} representation and formulates the free-boundary \glspl{pde}, namely obstacle-type \glspl{bvp} and the time-dependent \gls{sp}. Section~\ref{sec:methodology} details the proposed \gls{kan}-based physics-informed approach, emphasizing the residual-based loss and the training procedure. Section~\ref{sec:experimental_findings} reports numerical results for the \gls{eop}, \gls{plop}, and \gls{sp}, including comparisons with \gls{pinn} and \gls{resnet} baselines. Finally, Section~\ref{sec:conclusion} summarizes the work and discusses potential directions for future research. The Appendices~\ref{appendix:abbreviation_notations}-\ref{appendix:comparison_study_kan_parameters} contain additional implementation details, hyperparameter settings, and supplementary numerical comparisons.
\section{Problem Statement}
\label{sec:problem-statement}
\subsection{Kolmogorov-Arnold Network Architecture}
The Kolmogorov–Arnold representation theorem states that any multivariate continuous function on a bounded domain can be expressed as a finite composition of univariate continuous functions and summation operations \cite{arnol1957functions,kolmogorov1961representation,liu2024kan}.

\begin{figure}[!htbp]
\centering
\resizebox{\columnwidth}{!}{%
\begin{tikzpicture}[
    scale = 1.0,
    >=Latex,
    node/.style={circle, fill=black, inner sep=3.1pt},
    fun/.style={draw, rectangle, minimum width=0.65cm, minimum height=0.38cm, inner sep=1pt},
    edgefun/.style={draw=black, line width=0.45pt},
    post/.style={blue!70!black, font=\tiny},
    philab/.style={red!75!black, font=\tiny},
    nodelab/.style={font=\scriptsize},
    splinelab/.style={font=\tiny},
    every node/.style={font=\tiny}
]


\node[node,label=below:{\nodelab $a_{0,1}$}] (a01) at (0,0) {};
\node[node,label=below:{\nodelab $a_{0,2}$}] (a02) at (4.0,0) {};

\node[node,label=left:{\nodelab $a_{1,1}$}] (a11) at (-0.7,2.0) {};
\node[node,label=below:{\nodelab $a_{1,2}$}] (a12) at (0.8,2.0) {};
\node[node,label=above:{\nodelab $a_{1,3}$}] (a13) at (2.1,2.0) {};
\node[node,label=below:{\nodelab $a_{1,4}$}] (a14) at (3.4,2.0) {};
\node[node,label=below:{\nodelab $a_{1,5}$}] (a15) at (4.5,2.0) {};

\node[node,label=above:{\nodelab $a_{2,1}$}] (a21) at (2.1,4.4) {};

\foreach \j/\x/\y in {
1/-1.0/0.95,
2/-0.1/0.95,
3/0.8/0.95,
4/1.7/0.95,
5/2.6/0.95,
6/3.5/0.95,
7/4.4/0.95,
8/5.3/0.95}
{
    \node[fun] (f0\j) at (\x,\y) {};
}

\node[philab] at (f01) {$\phi_{0,1,1}$};
\node[philab] at (f02) {$\phi_{0,2,1}$};
\node[philab] at (f03) {$\phi_{0,3,1}$};
\node[philab] at (f04) {$\phi_{0,4,1}$};
\node[philab] at (f05) {$\phi_{0,5,1}$};
\node[philab] at (f06) {$\phi_{0,1,2}$};
\node[philab] at (f07) {$\phi_{0,3,2}$};
\node[philab] at (f08) {$\phi_{0,5,2}$};

\foreach \j/\x in {1/-0.7,2/0.8,3/2.1,4/3.4,5/4.9}
{
    \node[fun] (f1\j) at (\x,3.05) {};
    \node[philab] at (f1\j) {$\phi_{1,1,\j}$};
}

\draw[edgefun] (a01) -- (f01);
\draw[edgefun] (a01) -- (f02);
\draw[edgefun] (a01) -- (f03);
\draw[edgefun] (a01) -- (f04);
\draw[edgefun] (a01) -- (f05);

\draw[edgefun] (a02) -- (f06);
\draw[edgefun] (a02) -- (f07);
\draw[edgefun] (a02) -- (f08);

\draw[edgefun] (f01) -- (a11);
\draw[edgefun] (f02) -- (a12);
\draw[edgefun] (f03) -- (a13);
\draw[edgefun] (f04) -- (a14);
\draw[edgefun] (f05) -- (a15);

\draw[edgefun] (f06) -- (a11);
\draw[edgefun] (f06) -- (a12);
\draw[edgefun] (f07) -- (a13);
\draw[edgefun] (f07) -- (a14);
\draw[edgefun] (f08) -- (a15);

\foreach \j in {1,...,5}
{
    \draw[edgefun] (a1\j) -- (f1\j);
    \draw[edgefun] (f1\j) -- (a21);
}

\node[post] at (-0.25,2.18) {$\tilde{a}_{0,1,1}$};
\node[post] at (0.55,2.18) {$\tilde{a}_{0,2,1}$};
\node[post] at (1.55,2.18) {$\tilde{a}_{0,3,1}$};
\node[post] at (3.00,2.18) {$\tilde{a}_{0,4,2}$};
\node[post] at (4.35,2.18) {$\tilde{a}_{0,5,2}$};

\node[post] at (-0.4,3.65) {$\tilde{a}_{1,1,1}$};
\node[post] at (0.75,3.65) {$\tilde{a}_{1,1,2}$};
\node[post] at (2.05,3.65) {$\tilde{a}_{1,1,3}$};
\node[post] at (3.25,3.65) {$\tilde{a}_{1,1,4}$};
\node[post] at (4.55,3.65) {$\tilde{a}_{1,1,5}$};


\draw[dashed, thick] (8.95,2.55) ellipse (3.6cm and 2.6cm);

\draw[thick, ->]
    (f15.east) .. controls (5.95,3.55) and (6.35,4.35) .. (6.95,4.35);

\node at (6.80,4.15) {\scriptsize $k=3$};

\begin{scope}[shift={(7.20,3.30)}]
    \draw[->] (-1.7,0) -- (2.75,0) node[right] {$s$};

    \foreach \x/\lab in {
        -1.4/$s_{-3}$,
        -0.9/$s_{-2}$,
        -0.4/$s_{-1}$,
         0.1/$s_0$,
         0.65/$s_1$,
         1.2/$s_2$,
         1.75/$s_3$,
         2.3/$s_4$}
    {
        \draw (\x,0.05) -- (\x,-0.05) node[below] {\tiny \lab};
    }

    \foreach \s/\c in {
        -1.1/red!40,
        -0.6/orange!45,
        -0.1/green!45,
         0.4/blue!35,
         0.9/purple!35,
         1.4/cyan!40,
         1.9/red!30}
    {
        \draw[\c, smooth, domain=-0.45:0.45, samples=50, variable=\u]
        plot ({\s+\u},{0.45*exp(-8*\u*\u)});
    }

    \draw[very thick, smooth, domain=-0.6:1.85, samples=80, variable=\u]
        plot ({\u},{0.65 + 0.55*sin(150*\u) + 0.18*cos(330*\u)});

    \node[splinelab] at (1.35,1.55) {$\phi(s)$};
    \node[splinelab] at (2.55,0.98) {$\displaystyle \phi(s)=\sum_{i=0}^{7} c_i B_i(s)$};
    \node[splinelab] at (3.00,-0.30) {$G_1=5$};
\end{scope}

\draw[->, thick] (8.95,2.85) -- (8.95,2.20);
\node[right] at (9.03,2.65) {\tiny grid extension};

\begin{scope}[shift={(7.20,1.25)}]
    \draw[->] (-1.7,0) -- (2.75,0) node[right] {$s$};

    \foreach \x in {-1.4,-1.1,-0.8,-0.5,-0.2,0.1,0.4,0.7,1.0,1.3,1.6,1.9,2.2}
    {
        \draw (\x,0.05) -- (\x,-0.05);
    }

    \foreach \s/\c in {
        -1.35/red!40,
        -1.05/orange!45,
        -0.75/green!45,
        -0.45/blue!35,
        -0.15/purple!35,
         0.15/cyan!40,
         0.45/red!30,
         0.75/orange!35,
         1.05/green!35,
         1.35/blue!30,
         1.65/purple!30,
         1.95/cyan!35}
    {
        \draw[\c, smooth, domain=-0.25:0.25, samples=50, variable=\u]
        plot ({\s+\u},{0.42*exp(-18*\u*\u)});
    }

    \draw[very thick, smooth, domain=-0.6:1.85, samples=100, variable=\u]
        plot ({\u},{0.58 + 0.50*sin(150*\u) + 0.16*cos(330*\u)});

    \node[splinelab] at (2.55,0.95) {$\displaystyle \phi(s)=\sum_{i=0}^{12} c_i B_i(s)$};
    \node[splinelab] at (3.00,-0.30) {$G_2=10$};
\end{scope}
\end{tikzpicture}%
}
\caption{ \textbf{Kolmogorov-Arnold network architecture and spline-based grid refinement.}
Left: the \gls{kan} represents the solution using trainable univariate edge functions \(\phi_{l,j,i}\), where \(a_{l,i}\) denotes the scalar value of node \(i\) in layer \(l\), and \(\tilde{a}_{l,j,i}=\phi_{l,j,i}(a_{l,i})\) denotes the corresponding post-activation. Each node value \(a_{l+1,j}\) is obtained by summing incoming post-activations. Right: each edge function is parameterized by a univariate B-spline expansion in the scalar variable \(s\). Refining the spline grid from \(G_1=5\) to \(G_2=10\) increases the local approximation resolution.
}
\label{fig:kan-tikz}
\end{figure}
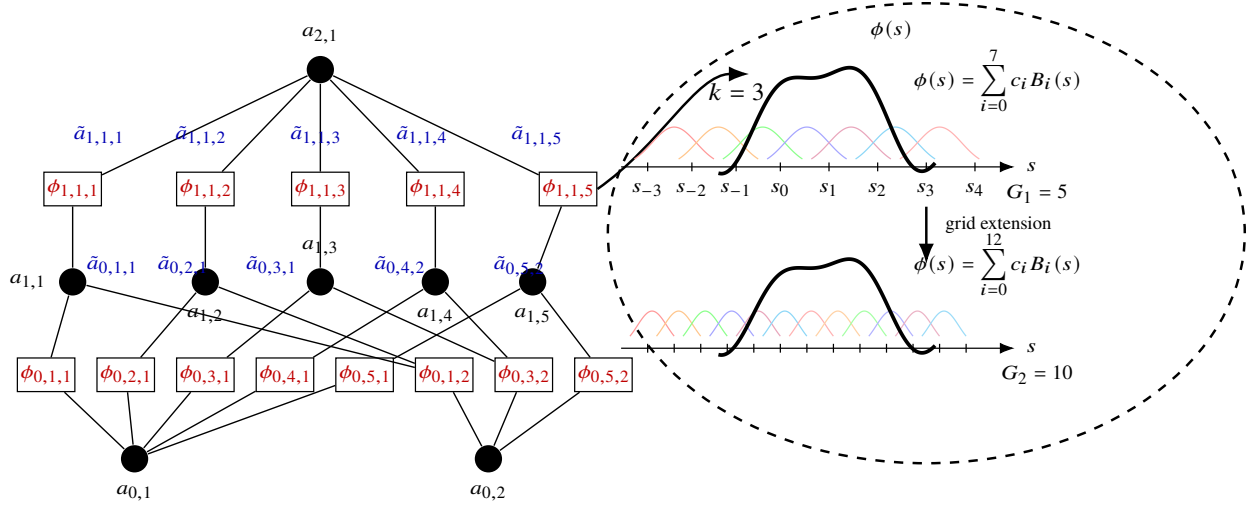

\begin{theorem}[Kolmogorov-Arnold Representation Theorem, see \cite{arnol1957functions, kolmogorov1961representation}]
Let \(h:[0,1]^d\to\mathbb{R}\) be a continuous multivariate function. There exist \(2d + 1\) continuous univariate functions \(\Phi^{(q)} : \mathbb{R} \to \mathbb{R}\) and \(d(2d+1)\) continuous univariate functions \(\phi^{(pq)} : [0,1] \to \mathbb{R}\) such that:
\begin{equation}
h (x_{1}, \dots, x_{d} )
=
\sum_{q=1}^{2d+1}
\Phi_q
\left(
\sum_{p=1}^{d}
\phi_{q,p}(x_p)
\right).
\label{eq:kan_approximation_formulation}
\end{equation}
\label{theorem:Arnold-Kolmogorov}
\end{theorem}
Motivated by Theorem~\ref{theorem:Arnold-Kolmogorov}, \glspl{kan} generalize the two-level Kolmogorov-Arnold decomposition to trainable networks with arbitrary depth and width to approximate multivariate functions by composing learnable one-dimensional functions placed on network edges. Specifically, we let the \gls{kan} architecture be specified by
\[
    [n_0,n_1,\ldots,n_L],
\]
where \(n_0\) is the input dimension, \(n_L=1\) is the output dimension, and \(L\) denotes the number of layers. For stationary problems, the input is \(\vec{x}\in \overline{\Omega}\subseteq\mathbb{R}^d\), where \(\Omega\) denotes the interior of the computational domain and \(\partial\Omega\) denotes its boundary, so \(\overline{\Omega}=\Omega\cup\partial\Omega\). For time-dependent problems, the input is \((\vec{x},t)\in\overline{\Omega}\times[0,T]\), so \(n_0=d+1\).
We denote by \(a_{l,i}\) the scalar value of node \(i\) in layer \(l\). In the input layer, we set
\[
    a_0 = \vec{x} \in  \overline{\Omega},
\]
for stationary problems and
\[
    a_0 = (\vec{x},t) \in \overline{\Omega} \times [0, T],
\]
for time-dependent problems. In a \gls{kan}, the connection from node \(i\) in layer \(l\) to node \(j\) in layer \(l+1\) is a learnable univariate function.
The next-layer node value is obtained by summing all edge functions:
\begin{equation}
    a_{l+1,j}
    =
    \sum_{i=1}^{n_l}
    \phi_{l,j,i}(a_{l,i}),
    \qquad j=1,\ldots,n_{l+1}.
    \label{eq:each_layer_node_summation_edge}
\end{equation}
Thus, each \gls{kan} layer can be written as
\[
    a_{l+1}=\Phi_l(a_l),
\]
where \(\Phi_l\) collects all learnable univariate edge functions between layers \(l\) and \(l+1\).
Each edge function is parameterized by a B-spline expansion,
\begin{equation}
    \phi_{l,j,i}(a_{l,i})
    =
    \sum_{r=0}^{G+k-1}
    c_{r}^{(l,j,i)} B_r(a_{l,i}),
    \label{eq:kan_bspline_edge}
\end{equation}
where \(B_r\) denotes the \(r\)-th univariate B-spline basis function, \(c_{r}^{(l,j,i)}\) is a trainable coefficient, \(G\) is the number of spline grid intervals, and \(k\) is the spline order.
Therefore, the unknown solution is approximated by a \gls{kan} denoted by \(u_{\theta}\). Figure~\ref{fig:kan-tikz} illustrates the \gls{kan} computational graph, in which every edge carries a learnable univariate function, and each node aggregates the incoming post-activations through summation. The right-hand panel shows the spline parameterization of an edge function and illustrates how grid refinement increases local resolution. For a stationary free-\gls{bvp} on \(\overline{\Omega}\subseteq\mathbb{R}^{d}\), we write
\begin{equation}
    u(\vec{x}) \approx u_{\theta}(\vec{x})
    =
    (\Phi_{L-1}\circ \Phi_{L-2}\circ \cdots \circ \Phi_{0})(\vec{x}).
    \label{eq:kan_solution_stationary_approximation}
\end{equation}
For a time-dependent free-\gls{bvp}, the input is augmented by time:
\begin{equation}
    u(\vec{x},t) \approx u_{\theta}(\vec{x},t)
    =
    (\Phi_{L-1}\circ \Phi_{L-2}\circ \cdots \circ \Phi_{0})(\vec{x},t),
    \label{eq:kan_solution_time_approximation}
\end{equation}
where \(\theta\) denotes all trainable spline coefficients of the \gls{kan}.
\subsection{Elliptic Obstacle Problem}
We first consider the \gls{eop}, which is a stationary free-\gls{bvp}, as in \cite{tran20151, ros2018obstacle, carstensen2015optimal}. Let \(\Omega\subset\mathbb{R}^{d}\) be a bounded domain with boundary \(\partial\Omega\). Given an obstacle function \(\psi\), a source term \(f\), and boundary data \(g\), the problem is to find \(u:\overline{\Omega}\to\mathbb{R}\) such that
\begin{equation}
\begin{cases}
    u(\vec{x}) \geq \psi(\vec{x}), & \vec{x}\in\Omega, \\
    -\Delta u(\vec{x}) \geq f(\vec{x}), & \vec{x}\in\Omega, \\
    \bigl(u(\vec{x})-\psi(\vec{x})\bigr)
    \bigl(-\Delta u(\vec{x})-f(\vec{x})\bigr) =0, & \vec{x}\in\Omega, \\
    u(\vec{s}) = g(\vec{s}), & \vec{s}\in\partial\Omega.
\end{cases}
\label{eq:obstacle-problem}
\end{equation}
The contact set is defined by
\begin{equation}
    \mathcal{C}=\{\vec{x}\in\Omega:u(\vec{x})=\psi(\vec{x})\}.
\end{equation}
The free boundary separates the contact and non-contact regions. In the non-contact region \(\Omega\setminus\mathcal{C}\), we have \(u>\psi\), so the complementarity condition implies \(-\Delta u=f\). In the contact set \(\mathcal{C}\), the solution satisfies \(u=\psi\). We next consider the \gls{plop}, which preserves this complementarity structure but replaces the linear Laplacian with a nonlinear diffusion operator \cite{pang2019fpinns}. For \(p>1\), the \(p\)-Laplacian operator is defined by
\begin{equation}
-\Delta_p u= -\nabla\cdot\left(|\nabla u|^{p-2}\nabla u\right).
\label{eq:p-laplacian-operator}
\end{equation}
This nonlinear formulation involves a diffusion term that depends on the gradient magnitude \(|\nabla u|\).
\subsection{Stefan Problem: Time-Dependent Free Boundary}
We next consider the one-phase \gls{sp}, a time-dependent free-\gls{bvp}. The temperature field is coupled with a moving interface. For \(d=2\), we represent the interface by the free-surface function \(s(x_2,t)\), and the moving domain is
\begin{equation}
\Omega(t)=\left\{(x_1,x_2)\in\Omega:0<x_1<s(x_2,t),\ 0<x_2<1 \right\}.
\end{equation}
The temperature field \(u(x_1,x_2,t)\) satisfies the heat equation in the moving domain for any \((x_1, x_2) \in \Omega(t)\) and \(t \in [0, T]\),
\begin{equation}
\partial_t u -\Delta u=0,
\label{eq:stefan_heat_equation}
\end{equation}
On the moving interface \(x_1=s(x_2,t)\), the phase-transition condition is
\begin{equation}
u(s(x_2,t),x_2,t)=0, \quad 0<x_2<1,\quad t\in(0,T],
\label{eq:stefan_interface_temperature}
\end{equation}
and the Stefan condition couples the normal velocity of the interface with the heat flux. For the graph \(x_1=s(x_2,t)\), the outward unit normal and normal velocity are
\[
\vec{n}
=
\frac{(1,-\partial_{x_2}s)}
{\sqrt{1+(\partial_{x_2}s)^2}},
\qquad
V_n
=
\frac{\partial_t s}
{\sqrt{1+(\partial_{x_2}s)^2}}.
\]
The one-phase Stefan condition \(V_n=-\nabla u\cdot\vec{n}\) is therefore equivalent to
\begin{equation}
\partial_t s(x_2,t)
+
\partial_{x_1}u(s(x_2,t),x_2,t)
-
\partial_{x_2}s(x_2,t)\,
\partial_{x_2}u(s(x_2,t),x_2,t)
=0,
\quad 0<x_2<1.
\label{eq:stefan_free_surface_condition}
\end{equation}
The problem is completed by appropriate initial and boundary conditions. In the proposed method, one \gls{kan} approximates the temperature field and a second \gls{kan} approximates the moving interface.

\section{Methodology}
\label{sec:methodology}
We use a residual-based physics-informed framework to train the proposed \gls{kan} solver. We represent the unknown solution by \(u_{\theta}\) and optimize the trainable parameters \(\theta\) using the Adam optimizer \cite{kingma2014adam}. The loss enforces the governing \gls{pde}, initial or boundary conditions, obstacle constraints, complementarity conditions, and, when applicable, free-boundary interface conditions. All spatial and temporal derivatives in the residual terms are computed using automatic differentiation \cite{paszke2017automatic,baydin2018automatic}.
\subsection{Residual-Based Loss Formulation}
We let \(u_{\theta}\) denote the \gls{kan} approximation as in Eq.~\eqref{eq:kan_solution_stationary_approximation}. Instead of enforcing Dirichlet boundary conditions through a penalty term, we impose them directly by constructing a solution that satisfies the boundary condition exactly. For a stationary free-\gls{bvp} with boundary function \(g\) on \(\partial\Omega\) as in Eq.~\eqref{eq:obstacle-problem}, we write
\begin{equation}
    u_{\theta}(\vec{x})
    =
    g(\vec{x})
    +
    \rho(\vec{x})\,N_{\theta}(\vec{x}),
    \qquad \vec{x}\in\overline{\Omega},
    \label{eq:hard_bc_stationary}
\end{equation}
where \(N_{\theta}\) is the \gls{kan} output and \(\rho:\overline{\Omega}\to\mathbb{R}\) is a boundary mask satisfying
\begin{equation}
    \rho(\vec{x})=0,
    \qquad \vec{x}\in\partial\Omega.
\end{equation}
Therefore, \(u_{\theta}(\vec{x})=g(\vec{x})\) holds exactly on \(\partial\Omega\).

For the \gls{eop}, we define the \gls{pde} residual
\begin{equation}
    R_{\mathrm{eop}}(\vec{x};\theta)
    =
    -\Delta u_{\theta}(\vec{x})-f(\vec{x}),
    \qquad \vec{x}\in\Omega.
    \label{eq:eop_residual}
\end{equation}
The inequality constraint \(-\Delta u\geq f\) is enforced through the residual violation
\begin{equation}
    \mathcal{L}_{\mathrm{pde}}
    =
    \frac{|\Omega|}{m_{\Omega}}
    \sum_{m=1}^{m_{\Omega}}
    \left[
    \max\{- R_{\mathrm{eop}}(\vec{x}_m;\theta),0\}
    \right]^2,
    \label{eq:eop_pde_loss}
\end{equation}
where \(\{\vec{x}_m\}_{m=1}^{m_{\Omega}}\subset\Omega\) are interior collocation points, and \(|\Omega|\) is the Lebesgue measure of the domain \(\Omega\). The obstacle constraint \(u\geq\psi\) is imposed by
\begin{equation}
    \mathcal{L}_{\mathrm{obs}}
    =
    \frac{|\Omega|}{m_{\Omega}}
    \sum_{m=1}^{m_{\Omega}}
    \left[
    \max\{\psi(\vec{x}_m)-u_{\theta}(\vec{x}_m),0\}
    \right]^2.
    \label{eq:eop_obstacle_loss}
\end{equation}
The complementarity condition is enforced by
\begin{equation}
    \mathcal{L}_{\mathrm{comp}}
    =
    \frac{|\Omega|}{m_{\Omega}}
    \sum_{m=1}^{m_{\Omega}}
    \left[
    \bigl(u_{\theta}(\vec{x}_m)-\psi(\vec{x}_m)\bigr)
    R_{\mathrm{eop}}(\vec{x}_m;\theta)
    \right]^2.
    \label{eq:eop_complementarity_loss}
\end{equation}
The total loss for the \gls{eop} is
\begin{equation}
    \mathcal{L}_{\mathrm{EOP}}(\theta)
    =
    \lambda_{\mathrm{pde}}\mathcal{L}_{\mathrm{pde}}
    +
    \lambda_{\mathrm{obs}}\mathcal{L}_{\mathrm{obs}}
    +
    \lambda_{\mathrm{comp}}\mathcal{L}_{\mathrm{comp}},
    \label{eq:eop_total_loss}
\end{equation}
where \(\lambda_{\mathrm{pde}}, \lambda_{\mathrm{obs}},\lambda_{\mathrm{comp}}>0\) are weighting parameters that balance the contributions of the \gls{pde} inequality loss, obstacle-constraint loss, and complementarity loss, respectively.
The trainable parameters for the \gls{eop} are obtained by solving
\begin{equation}
    \theta^{*}_{\mathrm{EOP}}
    =
    \arg\min_{\theta}
    \mathcal{L}_{\mathrm{EOP}}(\theta).
    \label{eq:eop_training_objective}
\end{equation}
For the nonlinear \gls{plop}, we replace the Laplacian residual with the \(p\)-Laplacian residual
\begin{equation}
    R_{p}(\vec{x};\theta)= - \nabla\cdot \left( |\nabla u_{\theta}(\vec{x})|^{p-2} \nabla u_{\theta}(\vec{x}) \right) - f(\vec{x}).
\label{eq:plop_residual}
\end{equation}
The corresponding total loss is
\begin{equation}
    \mathcal{L}_{\mathrm{PLOP}}(\theta)
    =
    \lambda_{\mathrm{pde},p}\mathcal{L}_{\mathrm{pde},p}
    +
    \lambda_{\mathrm{obs}}\mathcal{L}_{\mathrm{obs}}
    +
    \lambda_{\mathrm{comp},p}\mathcal{L}_{\mathrm{comp},p},
    \label{eq:plop_total_loss}
\end{equation}
where \(\mathcal{L}_{\mathrm{pde},p}\) and \(\mathcal{L}_{\mathrm{comp},p}\) are obtained from
\eqref{eq:eop_pde_loss} and \eqref{eq:eop_complementarity_loss} by replacing \(R_{\mathrm{eop}}\) with \(R_p\).

For the time-dependent \gls{sp}, we use one \gls{kan}, \(u_{\theta}(\vec{x},t)\), for the temperature field and a second \gls{kan}, \(\hat{s}_{\eta}(x_2,t)\), for the moving interface \(x_1=s(x_2,t)\). The predicted moving domain is
\[
    \Omega_{\eta}(t)
    =
    \{(x_1,x_2)\in\Omega:0<x_1<\hat{s}_{\eta}(x_2,t),\ 0<x_2<1\}.
\]
To impose the boundary condition, we use a hard-constrained representation of the form
\begin{equation}
    u_{\theta}(\vec{x},t)
    =
    A(\vec{x},t)
    +
    \rho(\vec{x},t)\,N_{\theta}(\vec{x},t),
    \label{eq:hard_bc_stefan}
\end{equation}
where \(A(\vec{x},t)\) is chosen to satisfy the prescribed fixed spatial boundary data and \(\rho(\vec{x},t)\) vanishes on those fixed boundary segments. The initial conditions are enforced separately through \(\mathcal{L}_{\mathrm{ic}}\). The heat-equation residual is then
\begin{equation}
    R_{\mathrm{heat}}(\vec{x},t;\theta)
    =
    \partial_t u_{\theta}(\vec{x},t)
    -
    \Delta u_{\theta}(\vec{x},t).
    \label{eq:heat_residual}
\end{equation}
The interior residual loss is
\begin{equation}
    \mathcal{L}_{\mathrm{heat}}
    =
    \frac{| \Omega_t|}{m_{\Omega_t}}
    \sum_{m=1}^{m_{\Omega_t}}
    \left[
    R_{\mathrm{heat}}(\vec{x}_m,t_m;\theta)
    \right]^2.
    \label{eq:heat_loss}
\end{equation}
The initial temperature field \(u(\vec{x},0)=u_0(\vec{x})\) is enforced by
\begin{equation}
    \mathcal{L}_{\mathrm{ic},u}
    =
    \frac{|\Omega_0|}{m_{\Omega_0}}
    \sum_{m=1}^{m_{\Omega_0}}
    \left[
    u_{\theta}(\vec{x}_m,0)-u_0(\vec{x}_m)
    \right]^2,
    \qquad \vec{x}_m\in\Omega_0 .
    \label{eq:initial_temperature_loss}
\end{equation}
Since the moving interface is also learned, we impose the initial free-boundary location
\(\hat{s}_{\eta}(x_2,0)=s_0(x_2)\) by
\begin{equation}
    \mathcal{L}_{\mathrm{ic},s}
    =
    \frac{|\Gamma_0|}{m_{\Gamma_0}}
    \sum_{m=1}^{m_{\Gamma_0}}
    \left[
    \hat{s}_{\eta}(x_{2,m},0)-s_0(x_{2,m})
    \right]^2 .
    \label{eq:initial_interface_loss}
\end{equation}
The initial-condition loss is
\begin{equation}
    \mathcal{L}_{\mathrm{ic}}
    =
    \mathcal{L}_{\mathrm{ic},u}
    +
    \mathcal{L}_{\mathrm{ic},s}.
    \label{eq:total_initial_condition_loss}
\end{equation}
The interface temperature condition is enforced on the learned moving boundary by
\begin{equation}
    \mathcal{L}_{\mathrm{int}}
    =
    \frac{|\Gamma|}{m_{\Gamma}}
    \sum_{m=1}^{m_{\Gamma}}
    \left[
    u_{\theta}(\hat{s}_{\eta}(x_{2,m},t_m),x_{2,m},t_m)
    \right]^2.
    \label{eq:interface_temperature_loss}
\end{equation}
Let
\[
    \vec{x}^{\,\Gamma}_m
    =
    \bigl(\hat{s}_{\eta}(x_{2,m},t_m),x_{2,m}\bigr)
\]
denote a collocation point on the learned moving interface. Using the graph form in Eq.~\eqref{eq:stefan_free_surface_condition}, define
\[
R_{\mathrm{fb}}(x_{2,m},t_m;\theta,\eta)
=
\partial_t \hat{s}_{\eta}(x_{2,m},t_m)
+
\partial_{x_1}u_{\theta}(\vec{x}^{\,\Gamma}_m,t_m)
-
\partial_{x_2}\hat{s}_{\eta}(x_{2,m},t_m)\,
\partial_{x_2}u_{\theta}(\vec{x}^{\,\Gamma}_m,t_m).
\]
The Stefan-condition loss is
\begin{equation}
    \mathcal{L}_{\mathrm{fb}}
    =
    \frac{|\Gamma|}{m_{\Gamma}}
    \sum_{m=1}^{m_{\Gamma}}
    \left[
    R_{\mathrm{fb}}(x_{2,m},t_m;\theta,\eta)
    \right]^2.
    \label{eq:stefan_condition_loss}
\end{equation}

The total loss for the \gls{sp} is
\begin{equation}
    \mathcal{L}_{\mathrm{SP}}(\theta,\eta)
    =
    \lambda_{\mathrm{heat}}\mathcal{L}_{\mathrm{heat}}
    +
    \lambda_{\mathrm{ic}}\mathcal{L}_{\mathrm{ic}}
    +
    \lambda_{\mathrm{int}}\mathcal{L}_{\mathrm{int}}
    +
    \lambda_{\mathrm{fb}}\mathcal{L}_{\mathrm{fb}},
    \label{eq:stefan_total_loss}
\end{equation}
where \(\lambda_{\mathrm{heat}},\lambda_{\mathrm{ic}},\lambda_{\mathrm{int}},\lambda_{\mathrm{fb}}>0\) are the weighting parameters for the heat-equation residual, initial conditions, interface temperature condition, and Stefan free-boundary condition, respectively.
The trainable parameters are obtained by solving
\begin{equation}
    (\theta^{*}_{\mathrm{SP}},\eta^{*}_{\mathrm{SP}})
    =
    \arg\min_{\theta,\eta}
    \mathcal{L}_{\mathrm{SP}}(\theta,\eta).
    \label{eq:theta_optimization_Stefan_Problem}
\end{equation}
\subsection{Algorithmic Frameworks}
We summarize the training procedure in two algorithms. Algorithm~\ref{alg:stationary_obstacle_solver} covers both stationary obstacle-type problems by selecting either \(R_{\mathrm{eop}}\) or \(R_p\). Algorithm~\ref{alg:stefan_solver} treats the coupled temperature/interface optimization for the time-dependent \gls{sp}.

\begin{algorithm}[!htbp]
\caption{\gls{kan}-based solver for stationary obstacle-type problems}
\label{alg:stationary_obstacle_solver}
\begin{algorithmic}
\REQUIRE Domain \(\overline{\Omega}\subseteq\mathbb{R}^d\), obstacle \(\psi\), source term \(f\), boundary data \(g\), interior collocation points \(\{\vec{x}_m\}_{m=1}^{m_\Omega}\), \gls{kan} \(N_\theta\), boundary mask \(\rho\), loss weights, learning rate \(\alpha\), problem type \(q\in\{\mathrm{EOP},\mathrm{PLOP}\}\), and \(\texttt{maxIter}\).
\STATE Initialize the trainable parameters \(\theta^{(0)}\).
\STATE Define the hard-boundary trial solution
\[
u_\theta(\vec{x})
=
g(\vec{x})+\rho(\vec{x})N_\theta(\vec{x}),
\qquad
\rho|_{\partial\Omega}=0.
\]
\FOR{\(k=0,\ldots,\texttt{maxIter}-1\)}
    \STATE Evaluate \(u_{\theta^{(k)}}(\vec{x}_m)\) and the derivatives required by the selected differential operator using automatic differentiation.
    \IF{\(q=\mathrm{EOP}\)}
        \STATE Set \(R^{(k)}=R_{\mathrm{eop}}(\cdot;\theta^{(k)})\) and evaluate the losses in Eqs.~\eqref{eq:eop_pde_loss}--\eqref{eq:eop_complementarity_loss}.
        \STATE Form \(\mathcal{L}^{(k)}=\mathcal{L}_{\mathrm{EOP}}(\theta^{(k)})\) using Eq.~\eqref{eq:eop_total_loss}.
    \ELSE
        \STATE Set \(R^{(k)}=R_p(\cdot;\theta^{(k)})\), replace \(R_{\mathrm{eop}}\) by \(R_p\) in the inequality and complementarity losses, and form \(\mathcal{L}^{(k)}=\mathcal{L}_{\mathrm{PLOP}}(\theta^{(k)})\) using Eq.~\eqref{eq:plop_total_loss}.
    \ENDIF
    \STATE Update the parameters with Adam:
    \[
    \theta^{(k+1)}
    =
    \mathrm{Adam}\!\left(
    \theta^{(k)},
    \nabla_\theta \mathcal{L}^{(k)},
    \alpha
    \right).
    \]
\ENDFOR
\ENSURE \(u_{\theta^{(\texttt{maxIter})}}\).
\end{algorithmic}
\end{algorithm}

For the time-dependent \gls{sp}, the current interface determines the moving domain on which the heat-equation residual is evaluated. The temperature and interface networks are therefore updated jointly.

\begin{algorithm}[!htbp]
\caption{\gls{kan}-based solver for the time-dependent \gls{sp}}
\label{alg:stefan_solver}
\begin{algorithmic}
\REQUIRE Initial data \(u_0,s_0\), fixed boundary data, candidate space-time points, initial/interface collocation points, \gls{kan} models \(u_\theta\) and \(\hat{s}_\eta\), loss weights \(\lambda_{\mathrm{heat}},\lambda_{\mathrm{ic}},\lambda_{\mathrm{int}},\lambda_{\mathrm{fb}}\), learning rate \(\alpha\), and \(\texttt{maxIter}\).
\STATE Initialize \(\theta^{(0)}\) and \(\eta^{(0)}\).
\STATE Construct \(u_\theta=A+\rho N_\theta\) so that the prescribed fixed spatial boundary data are satisfied exactly.
\FOR{\(k=0,\ldots,\texttt{maxIter}-1\)}
    \STATE Define the current predicted domain
    \[
    \Omega_{\eta^{(k)}}(t)
    =
    \left\{
    (x_1,x_2)\in\Omega:
    0<x_1<\hat{s}_{\eta^{(k)}}(x_2,t),\
    0<x_2<1
    \right\}.
    \]
    \STATE Select/evaluate interior collocation points lying in \(\Omega_{\eta^{(k)}}(t)\), together with initial and interface points.
    \STATE Evaluate \(u_{\theta^{(k)}}\), \(\hat{s}_{\eta^{(k)}}\), and all required derivatives by automatic differentiation.
    \STATE Compute \(\mathcal{L}_{\mathrm{heat}}\), \(\mathcal{L}_{\mathrm{ic}}\), \(\mathcal{L}_{\mathrm{int}}\), and \(\mathcal{L}_{\mathrm{fb}}\) using Eqs.~\eqref{eq:heat_loss}, \eqref{eq:total_initial_condition_loss}, \eqref{eq:interface_temperature_loss}, and \eqref{eq:stefan_condition_loss}.
    \STATE Form
    \[
    \mathcal{L}^{(k)}_{\mathrm{SP}}
    =
    \lambda_{\mathrm{heat}}\mathcal{L}_{\mathrm{heat}}
    +
    \lambda_{\mathrm{ic}}\mathcal{L}_{\mathrm{ic}}
    +
    \lambda_{\mathrm{int}}\mathcal{L}_{\mathrm{int}}
    +
    \lambda_{\mathrm{fb}}\mathcal{L}_{\mathrm{fb}}.
    \]
    \STATE Jointly update both parameter sets:
    \[
    (\theta^{(k+1)},\eta^{(k+1)})
    =
    \mathrm{Adam}\!\left(
    (\theta^{(k)},\eta^{(k)}),
    \nabla_{\theta,\eta}\mathcal{L}^{(k)}_{\mathrm{SP}},
    \alpha
    \right).
    \]
\ENDFOR
\ENSURE \(u_{\theta^{(\texttt{maxIter})}}(\vec{x},t)\) and \(\hat{s}_{\eta^{(\texttt{maxIter})}}(x_2,t)\).
\end{algorithmic}
\end{algorithm}
\FloatBarrier

\section{Numerical Experiments}
\label{sec:experimental_findings}
We evaluate the proposed \gls{kan}-based physics-informed solver on the \gls{eop}, the nonlinear \gls{plop}, and the time-dependent \gls{sp}. The two stationary obstacle problems are solved using Algorithm~\ref{alg:stationary_obstacle_solver}, while the \gls{sp} is solved using Algorithm~\ref{alg:stefan_solver}. We evaluate solution accuracy, convergence of the residual losses, and the ability of the models to resolve contact regions and moving free boundaries. We report the relative \(L^2\) and \(L^\infty\) errors:
\begin{equation}
    E_{L^2}
    =
    \frac{\|u_{\theta}-u_{\mathrm{exact}}\|_{L^2(\overline{\Omega})}}
    {\|u_{\mathrm{exact}}\|_{L^2(\overline{\Omega})}},
    \label{eq:relative_error_l2_metrics}
\end{equation}
\begin{equation}
    E_{L^\infty}
    =
    \frac{\|u_{\theta}-u_{\mathrm{exact}}\|_{L^\infty(\overline{\Omega})}}
    {\|u_{\mathrm{exact}}\|_{L^\infty(\overline{\Omega})}}.
    \label{eq:relative_error_metric_inf}
\end{equation}
For time-dependent problems, the norms are evaluated over the corresponding space-time domain or at selected time snapshots. We also compare \gls{kan} models with \gls{pinn} and \gls{resnet} baselines to evaluate accuracy and parameter efficiency. Appendix~\ref{appendix:settings-pde} lists the problem settings, exact solutions, sampling points, network architectures, and optimization hyperparameters. Appendix~\ref{appendix:comparison_study_kan_parameters} presents an ablation study of spline grid size and network width.
\subsection{Radially Symmetric Obstacle Problem}\label{exp:radially-symmetric-obstacle}
For the \gls{eop}, we consider the computational domain
\[
    \overline{\Omega}=[0,1]^2.
\]
The obstacle problem is
\begin{equation}
\begin{cases}
    u(\vec{x}) \geq \psi(\vec{x}), & \vec{x}\in\Omega,\\
    -\Delta u(\vec{x}) \geq 0, & \vec{x}\in\Omega,\\
    \bigl(u(\vec{x})-\psi(\vec{x})\bigr)\bigl(-\Delta u(\vec{x})\bigr)=0,
    & \vec{x}\in\Omega,\\
    u(\vec{s})=g(\vec{s}), & \vec{s}\in\partial\Omega,
\end{cases}
\label{eq:appendix_eop_pde}
\end{equation}
where the source term is \(f(\vec{x})=0\), so the governing equation in the non-contact region is \(-\Delta u=0\).

\begin{figure}[!htbp]
    \centering
    \begin{subfigure}[t]{0.32\textwidth}
    \centering
    \includegraphics[width =\linewidth]{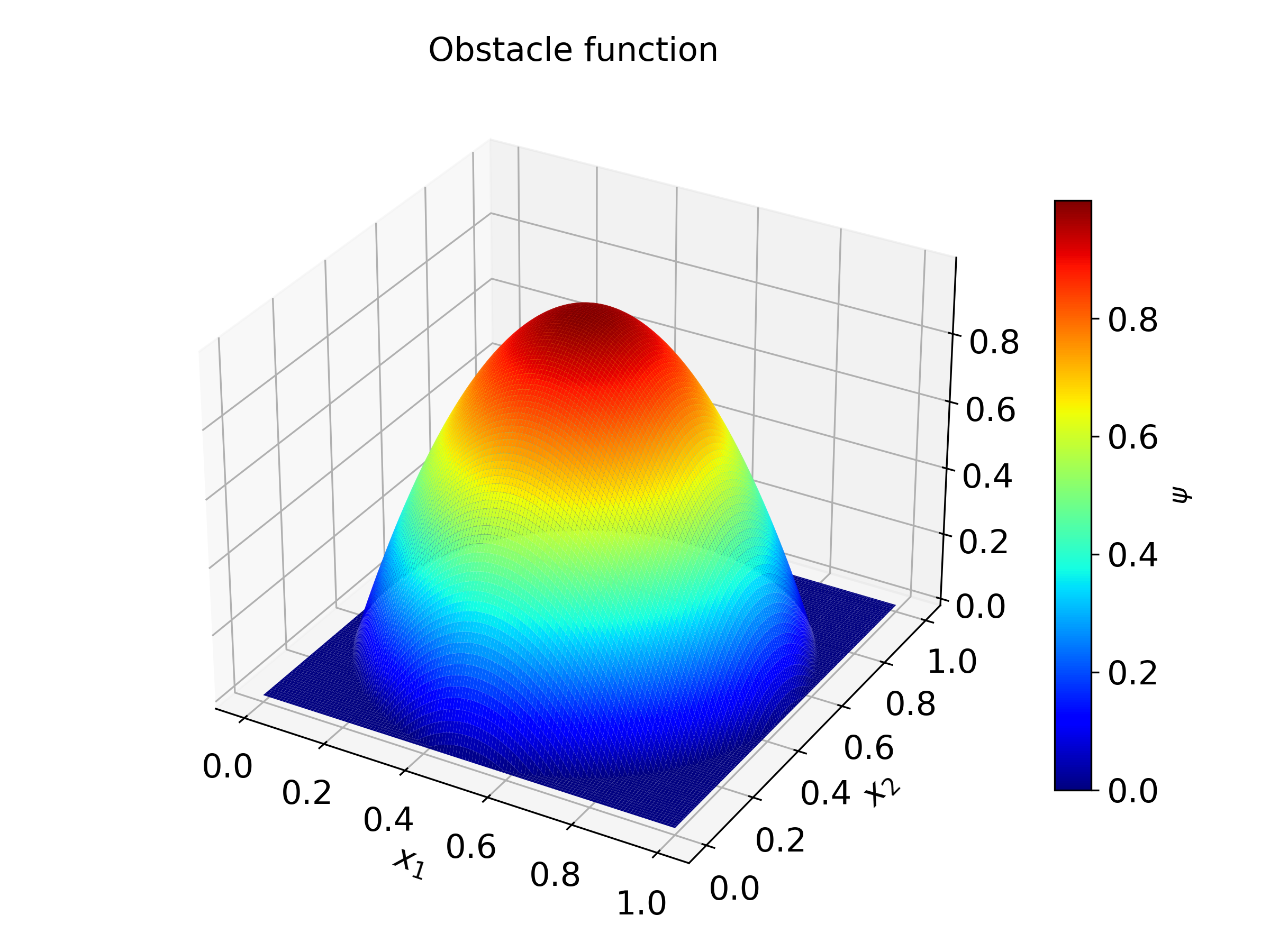}
    \caption{Obstacle function}
    \end{subfigure}
    \begin{subfigure}[t]{0.32\textwidth}
    \centering
    \includegraphics[width =\linewidth]{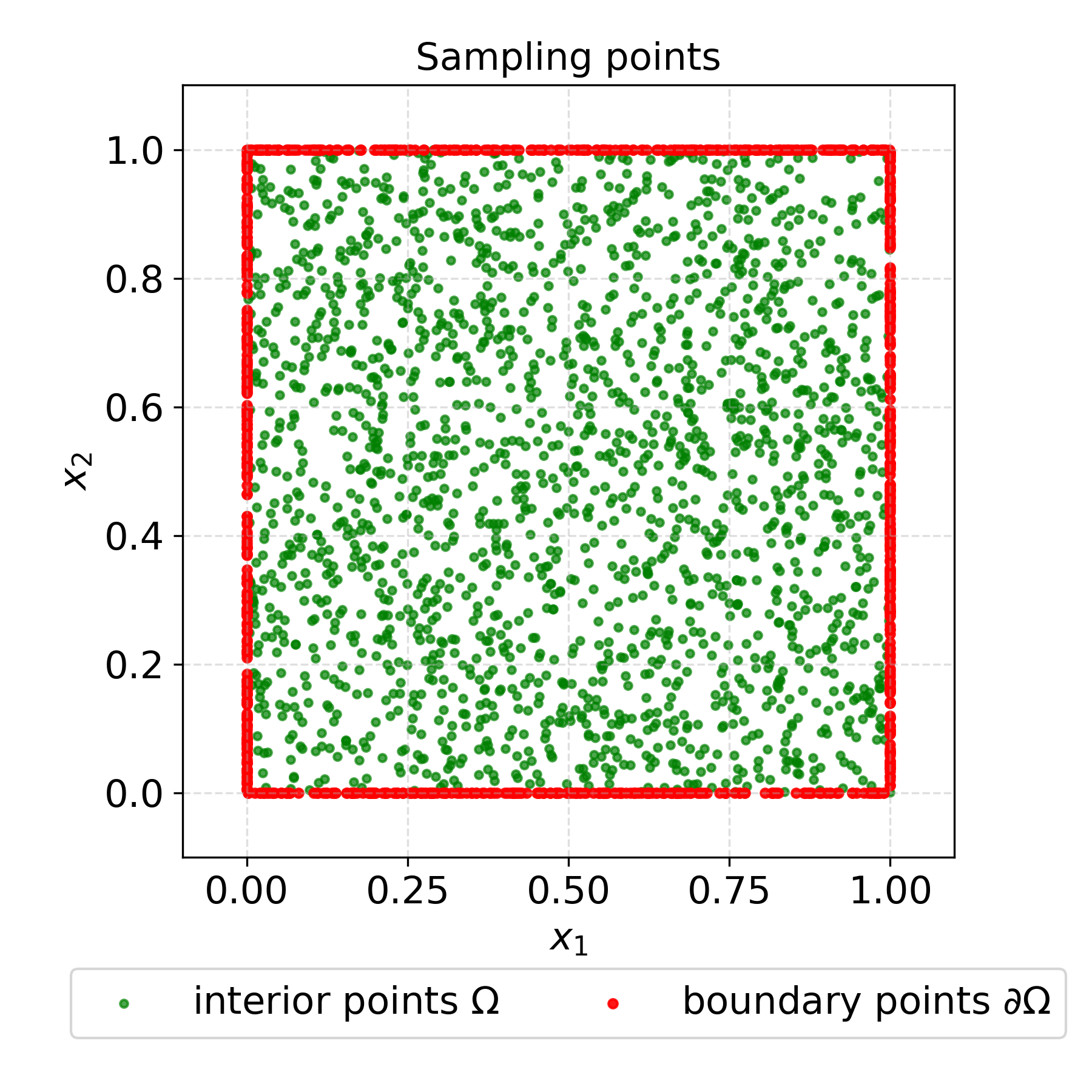}
    \caption{Sampled collocation points}\hfill
    \end{subfigure}
    \caption{\textbf{Obstacle function and collocation points.}
    Panel (a) shows the radial obstacle function \(\psi\) on \(\Omega=[0,1]^2\). Panel (b) shows the sampled collocation points, with interior points in green and boundary points on \(\partial\Omega\) in red.}
    \label{fig:obstacle_func_eop}
\end{figure}

We use a radially symmetric dome obstacle centered at \((1/2,1/2)\). Let
\[
    r(\vec{x})
    =
    \sqrt{(x_1-1/2)^2+(x_2-1/2)^2},
    \qquad
    r_c=\frac{1}{2}.
\]
The obstacle function is defined by
\begin{equation}
    \psi(\vec{x})
    =
    \begin{cases}
    1-\dfrac{r(\vec{x})^2}{r_c^2},
    & r(\vec{x})\leq r_c,\\[0.7em]
    0,
    & r(\vec{x})>r_c.
    \end{cases}
    \label{eq:appendix_eop_obstacle}
\end{equation}
The corresponding exact solution is
\begin{equation}
    u_{\mathrm{exact}}(\vec{x})
    =
    \begin{cases}
    1-\dfrac{r(\vec{x})^2}{r_c^2},
    & r(\vec{x})\leq r_*,\\[0.9em]
    1-\dfrac{r_*^2}{r_c^2}
    \left(
    1+2\log\left(\dfrac{r(\vec{x})}{r_*}\right)
    \right),
    & r(\vec{x})>r_*,
    \end{cases}
    \label{eq:appendix_eop_exact_solution}
\end{equation}
where \(r_*\approx0.260197\) denotes the contact radius. The Dirichlet boundary data are prescribed from the exact solution in Eq.~\eqref{eq:appendix_eop_exact_solution}:
\[
    g(\vec{s})=u_{\mathrm{exact}}(\vec{s}),
    \qquad \vec{s}\in\partial\Omega.
\]
Figure~\ref{fig:obstacle_func_eop} shows the obstacle function and the sampled collocation points.
\FloatBarrier
\begin{figure}[p]
    \centering
    \begin{subfigure}[t]{0.32\textwidth}
    \centering
    \includegraphics[width = \linewidth]{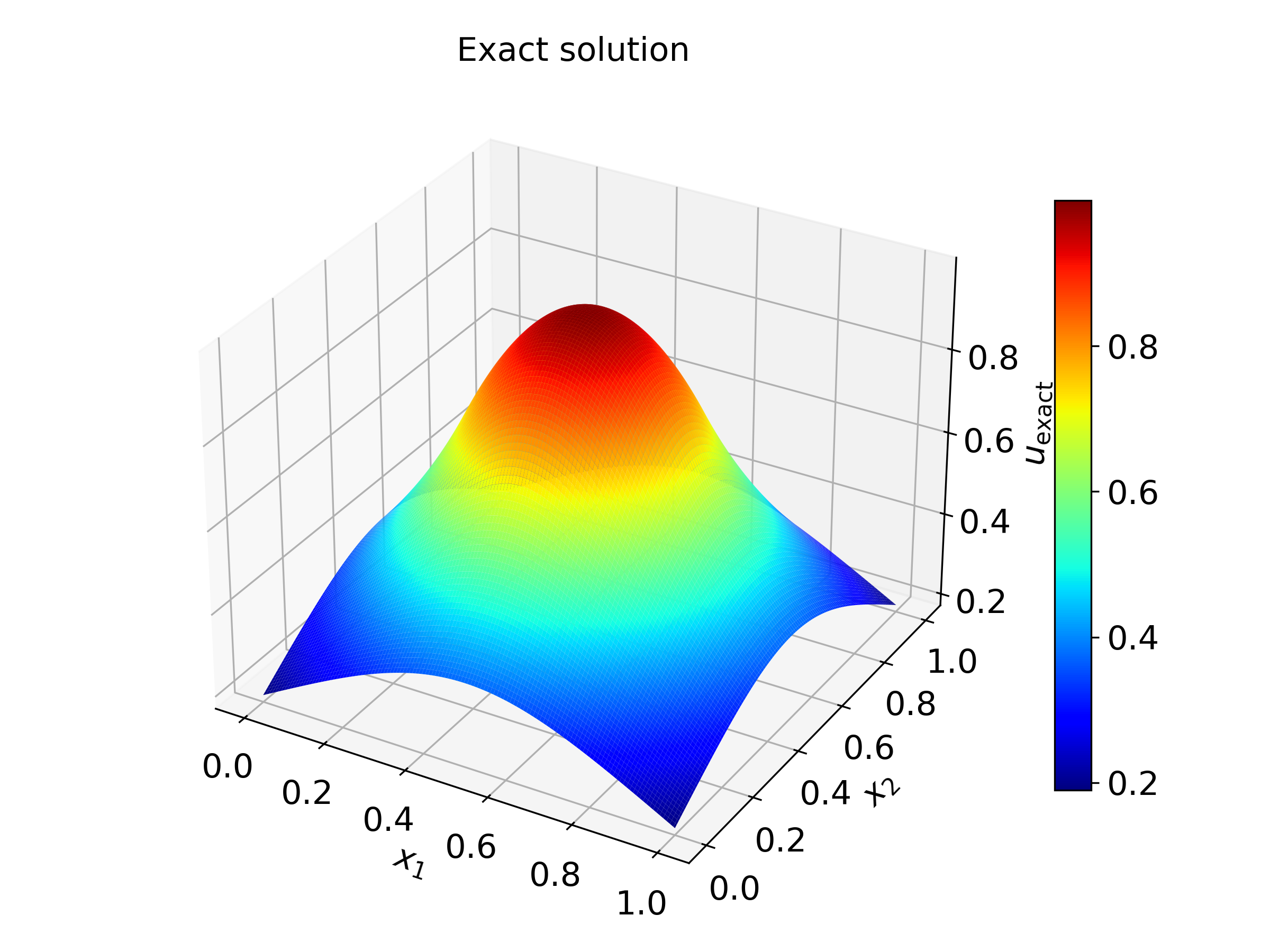}
    \caption{Exact solution}
    \end{subfigure}
    \begin{subfigure}[t]{0.32\textwidth}
    \centering
    \includegraphics[width = \linewidth]{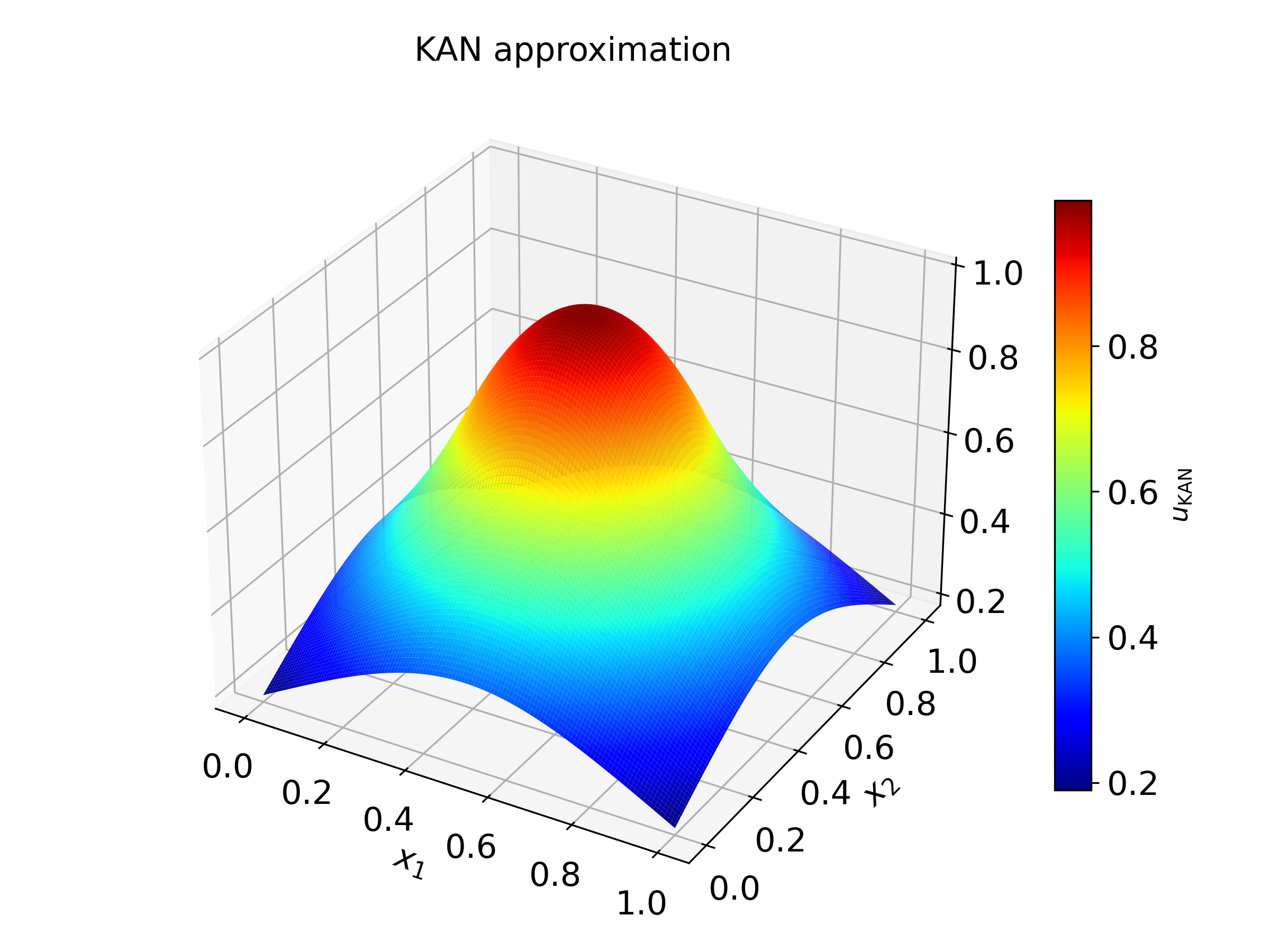}
    \caption{\(u_{\mathrm{KAN}}\) approximation}
    \end{subfigure}
    \begin{subfigure}[t]{0.32\textwidth}
    \centering
    \includegraphics[width = \linewidth]{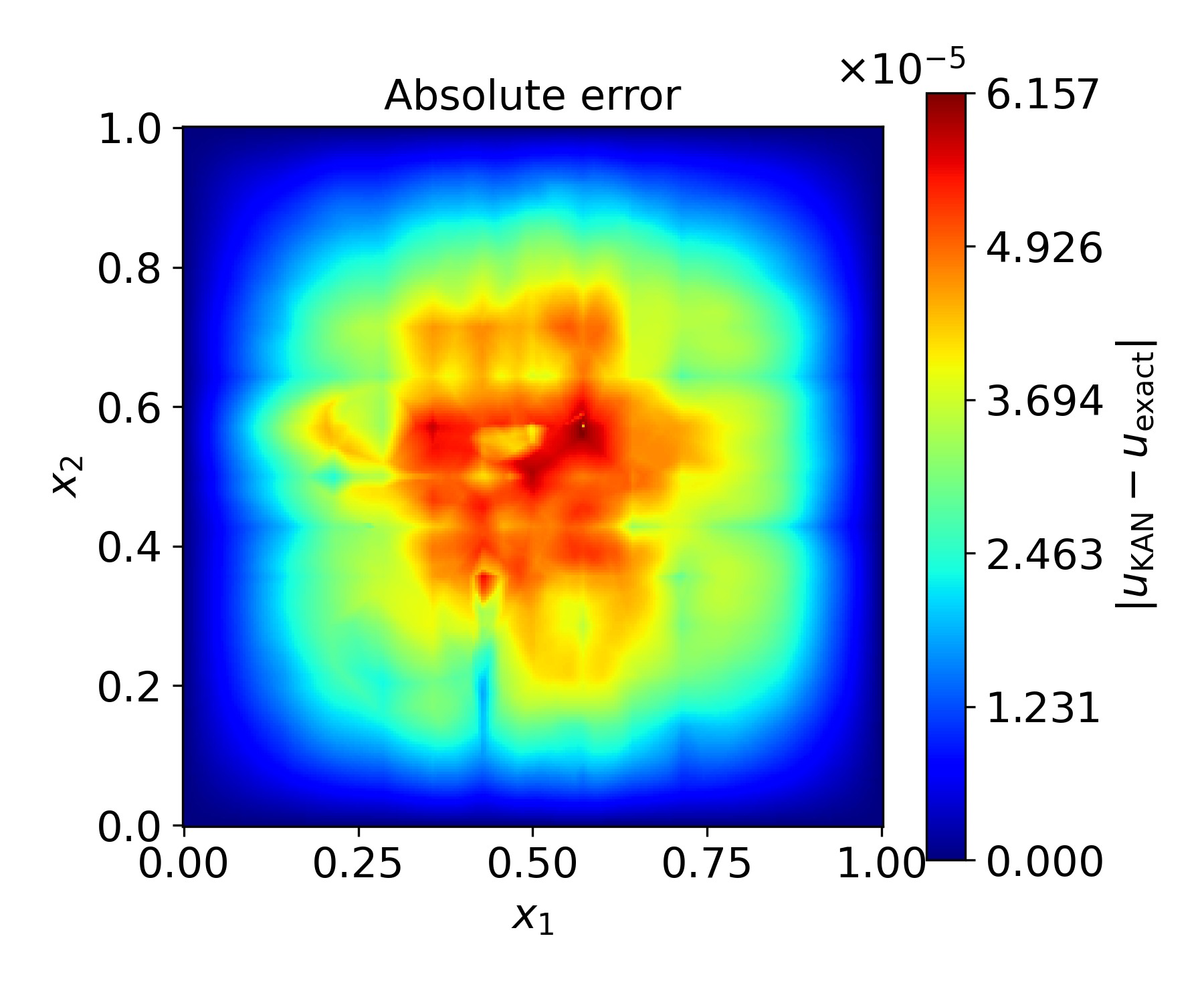}
    \caption{Error distribution}
    \end{subfigure}\hfill
    \vspace{0.2em}
    \begin{subfigure}[t]{0.5\textwidth}
    \centering
    \includegraphics[width = \linewidth]{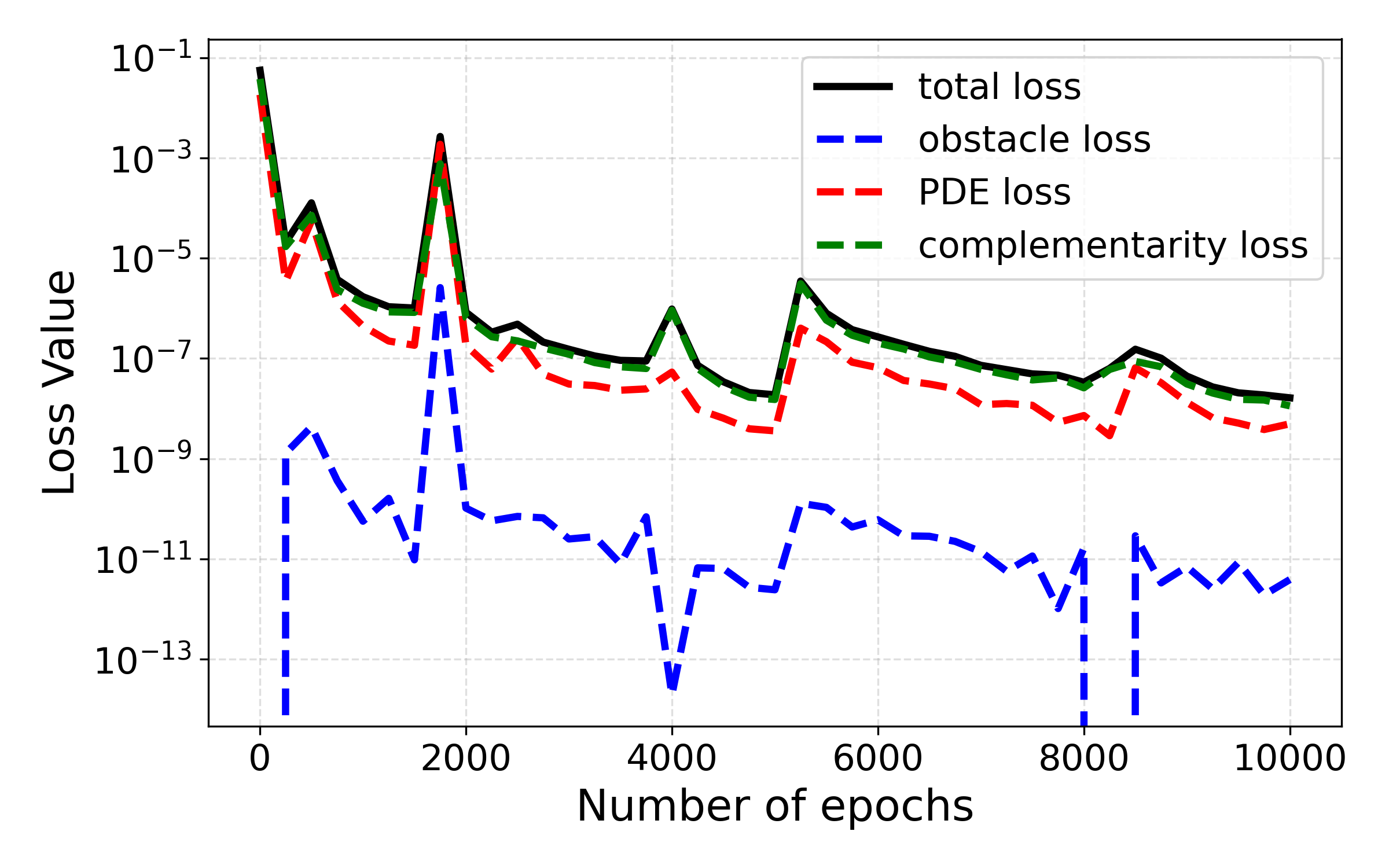}
    \caption{Loss components vs.\ epochs}
    \end{subfigure}\hfill
    \begin{subfigure}[t]{0.5\textwidth}
    \centering
    \includegraphics[width = \linewidth]{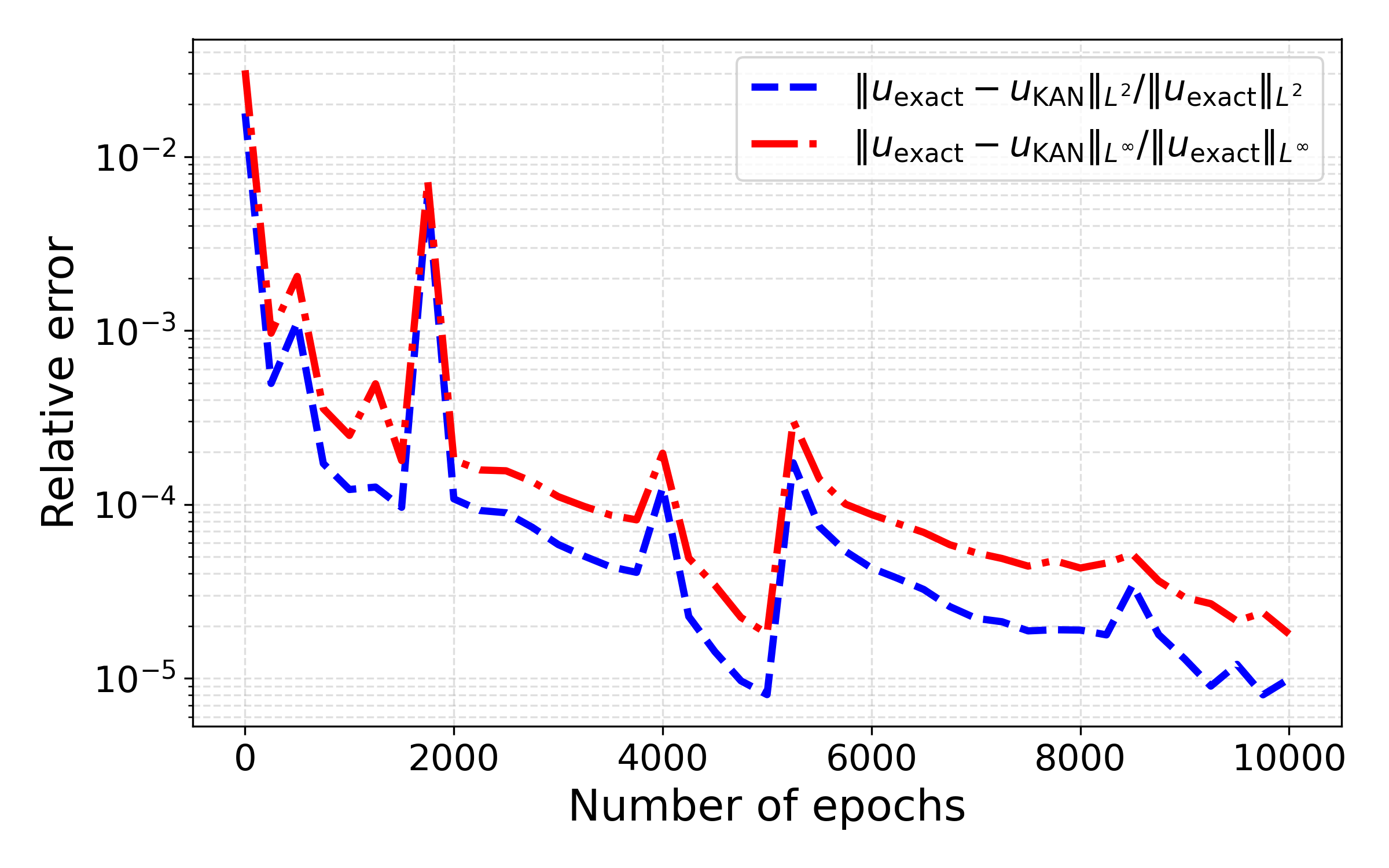}
    \caption{Relative error vs.\ epochs}
    \end{subfigure}\hfill
    \caption{\textbf{KAN approximation for the two-dimensional elastic obstacle problem.}
Comparison between the exact solution and the proposed \gls{kan} approximation on \(\Omega=[0,1]^2\). Panels (a) and (b) show the exact solution \(u_{\mathrm{exact}}\) and the learned approximation \(u_{\mathrm{KAN}}\), respectively. Panel (c) shows the absolute pointwise error \(|u_{\mathrm{KAN}}-u_{\mathrm{exact}}|\), together with the circular free-boundary interface of radius \(r^{*}=0.260\). Panel (d) reports the evolution of the total loss and its obstacle, PDE, and complementarity components during training. Panel (e) shows the relative \(L^{2}\) and \(L^{\infty}\) errors over the training epochs. The final relative \(L^2\) and \(L^{\infty}\) errors are \(4.858\times10^{-5}\) and \(6.517\times10^{-5}\), respectively.}
    \label{fig:kan_approximation_eop}
\end{figure}
\FloatBarrier
\begin{figure}[!htbp]
    \centering
    \includegraphics[width=\linewidth]{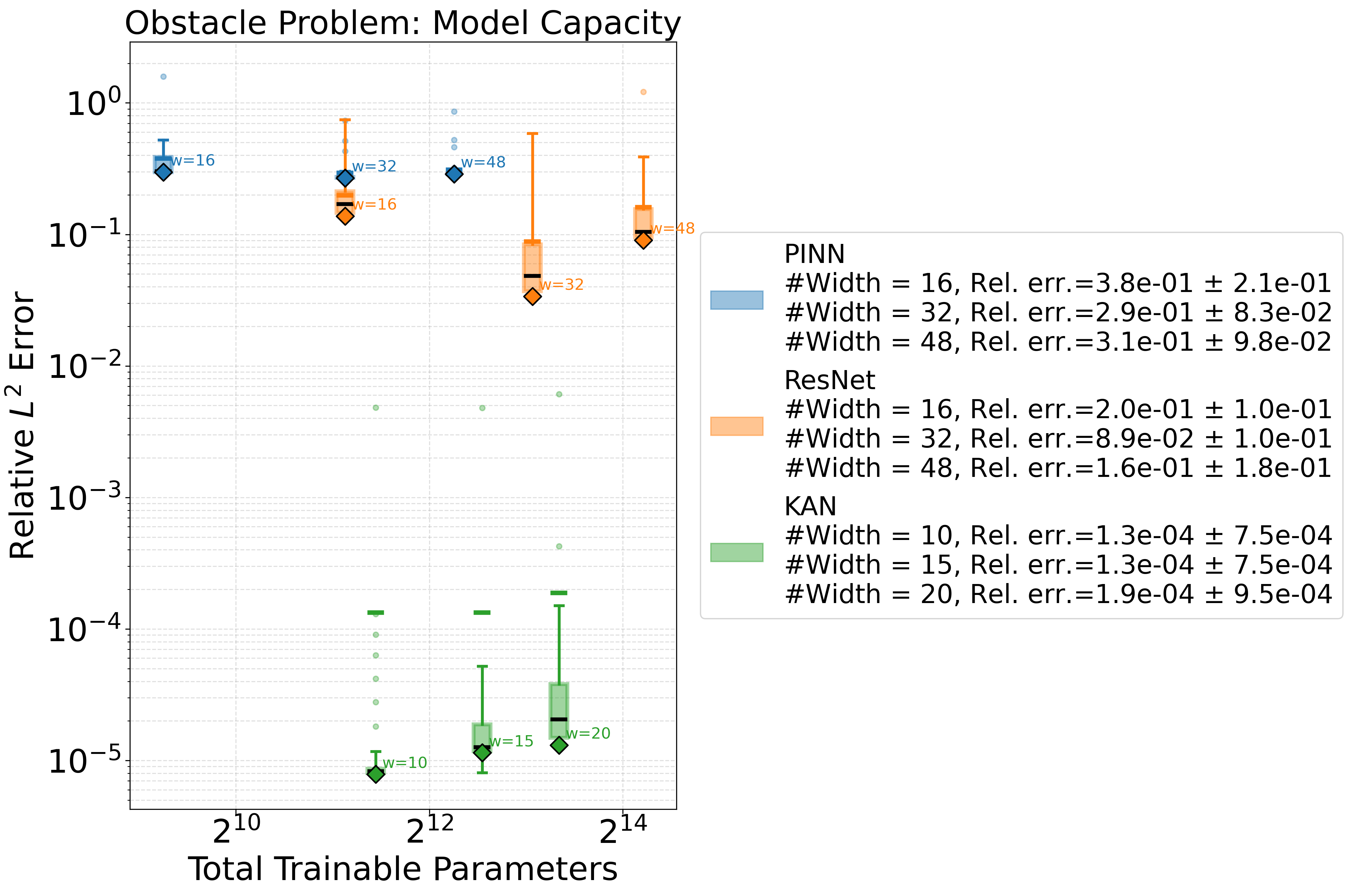}
    \caption{
    \textbf{Model-capacity comparison for the \gls{eop}.}
    Relative \(L^2\) errors are plotted against the total number of trainable parameters for \gls{pinn}, \gls{resnet}, and \gls{kan} models. Across the different architectures, \gls{kan} achieves lower errors with fewer parameters, demonstrating improved parameter efficiency for the obstacle problem.
}
\label{fig_comparison_methods_obstacle_problem}
\end{figure}

Figure~\ref{fig:kan_approximation_eop} shows close agreement between the \gls{kan} approximation and the exact solution. The model accurately resolves the circular contact region, and the pointwise error remains small across the computational domain. The final relative \(L^2\) and \(L^\infty\) errors are \(4.858\times10^{-5}\) and \(6.517\times10^{-5}\), respectively. Figure~\ref{fig_comparison_methods_obstacle_problem} further shows that the \gls{kan} models achieve lower relative \(L^2\) errors with fewer trainable parameters than the \gls{pinn} and \gls{resnet} baselines. The ablation study in Figure~\ref{fig_params_gridsize_study}(a) shows that the error depends jointly on spline grid size and network width. Increasing either quantity does not guarantee lower error.

\FloatBarrier
\subsection{p-Laplacian Obstacle Problem}\label{exp:p-Laplacian-Obstacle-Problem}
\begin{figure}[p]
    \centering
    \begin{subfigure}[t]{0.32\textwidth}
    \centering
    \includegraphics[width = \linewidth]{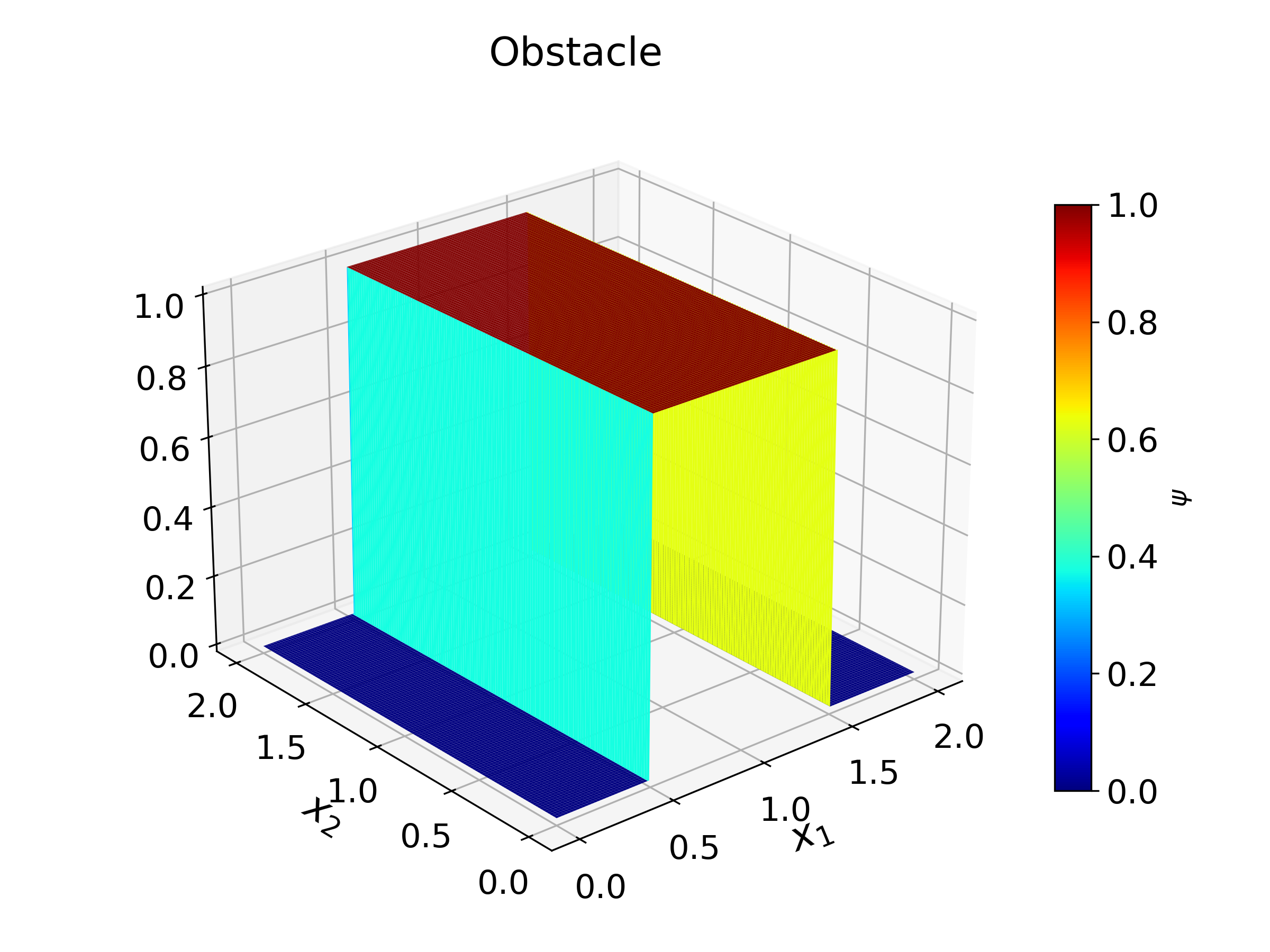}
    \caption{Obstacle function}
    \end{subfigure}
    \centering
    \begin{subfigure}[t]{0.32\textwidth}
    \centering
    \includegraphics[width = \linewidth]{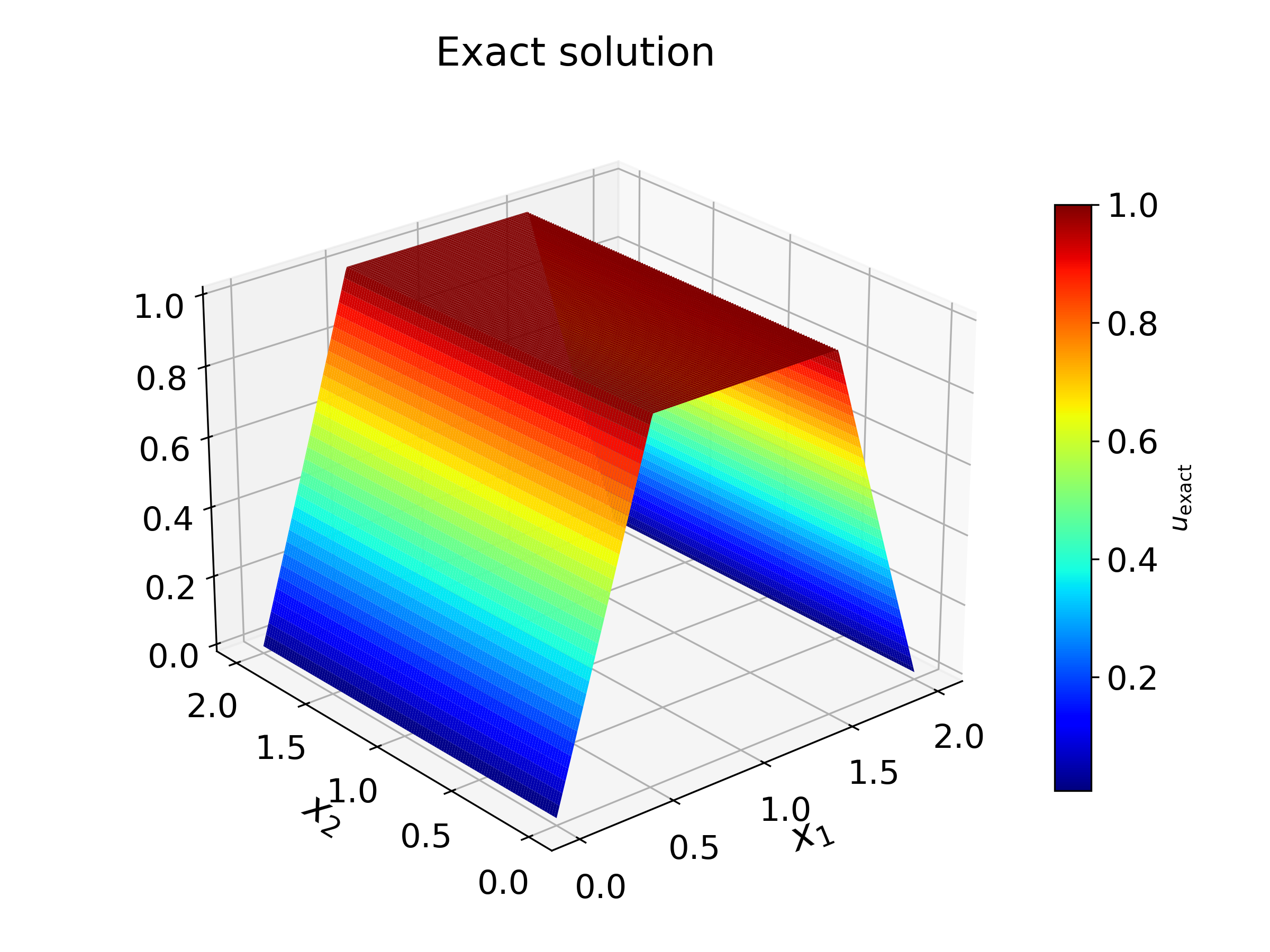}
    \caption{Exact solution}
    \end{subfigure}
    \centering
    \begin{subfigure}[t]{0.32\textwidth}
    \centering
    \includegraphics[width = \linewidth]{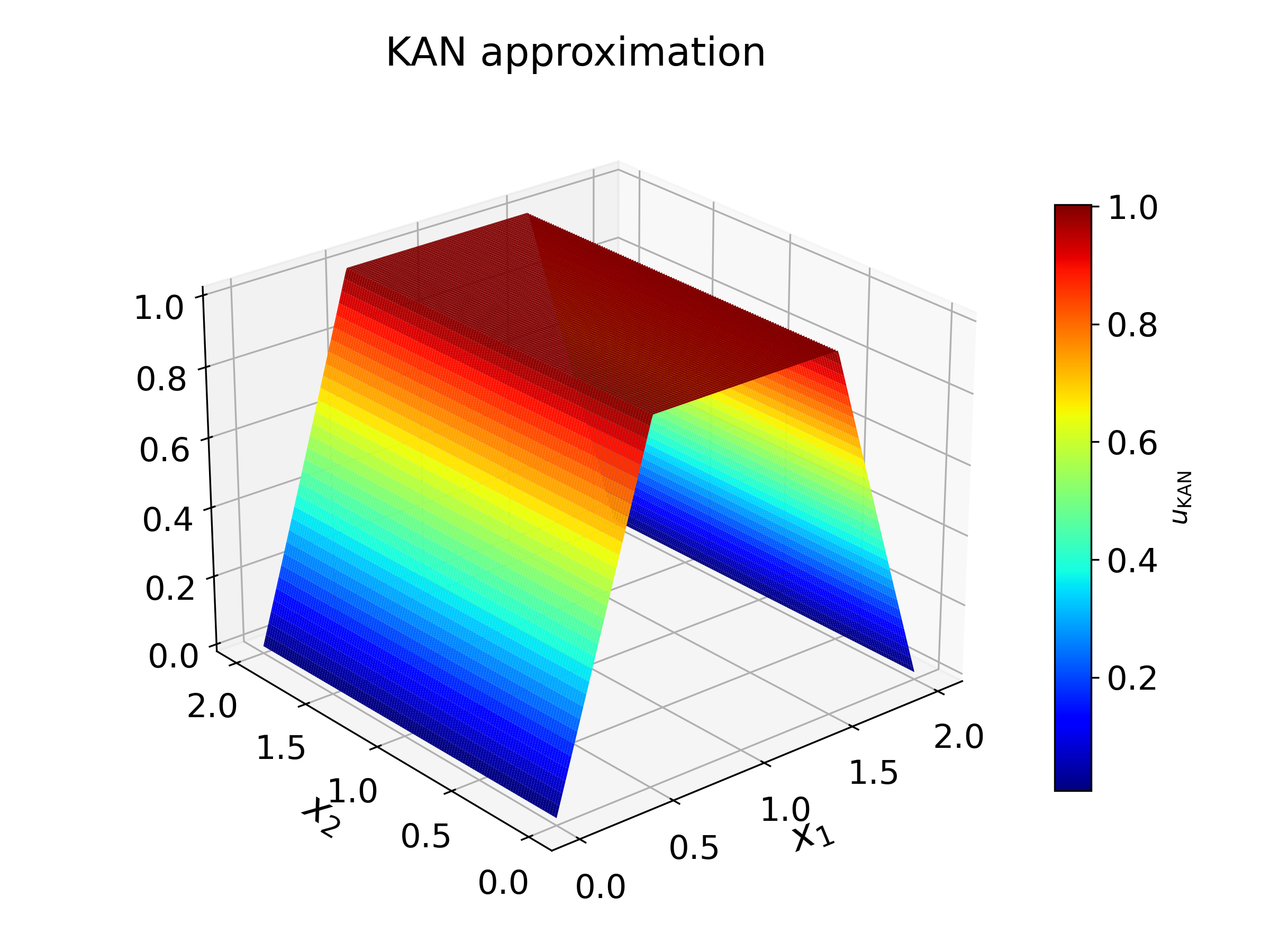}
    \caption{\gls{kan} approximation}
    \end{subfigure}\hfill
    \vspace{0.2em}
    \begin{subfigure}[t]{0.32\textwidth}
    \centering
    \includegraphics[width = \linewidth]{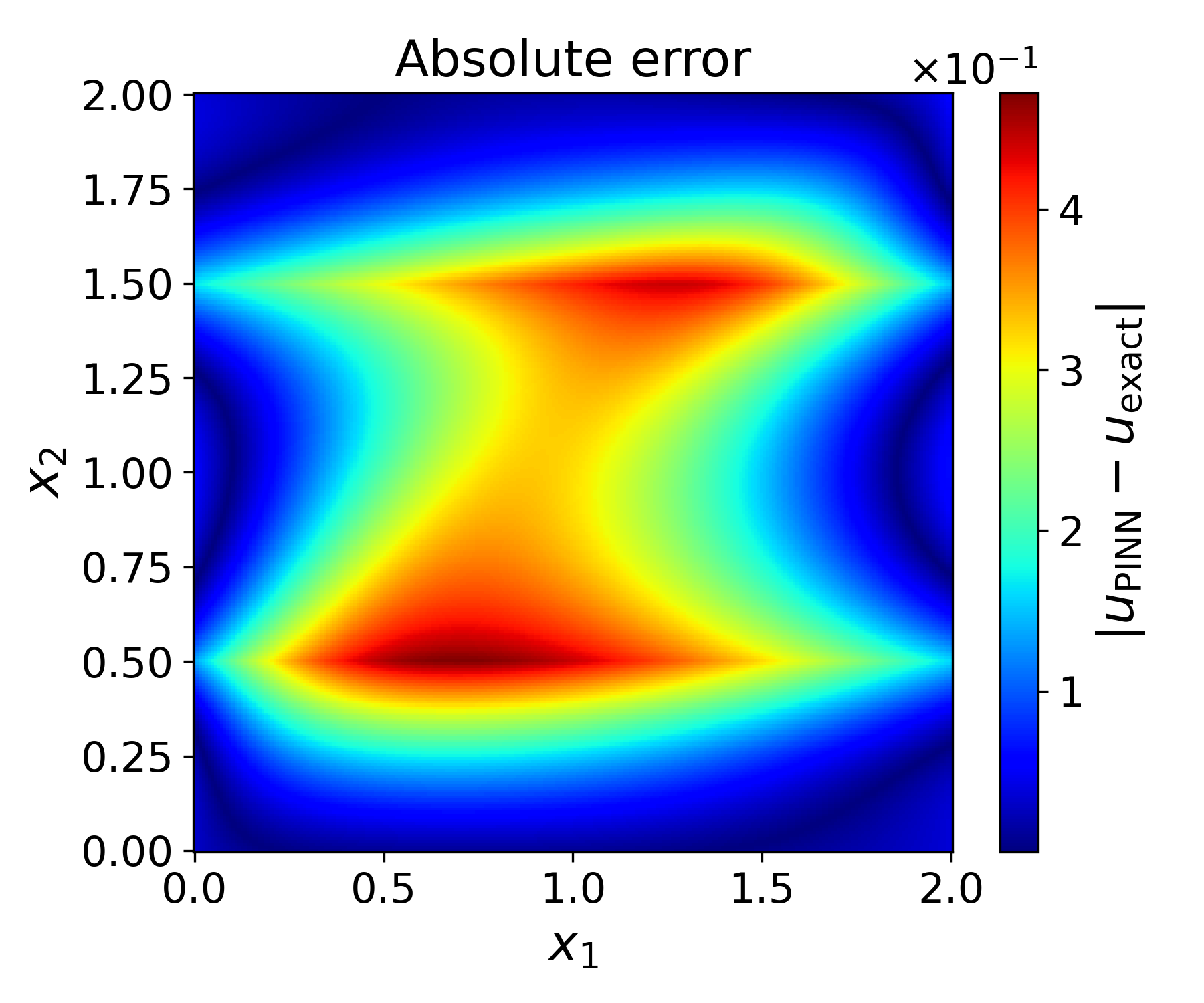}
    \caption{\(|u_{exact} - u_{\mathrm{PINN}}|\)}
    \end{subfigure}
    \begin{subfigure}[t]{0.32\textwidth}
    \centering
    \includegraphics[width = \linewidth]{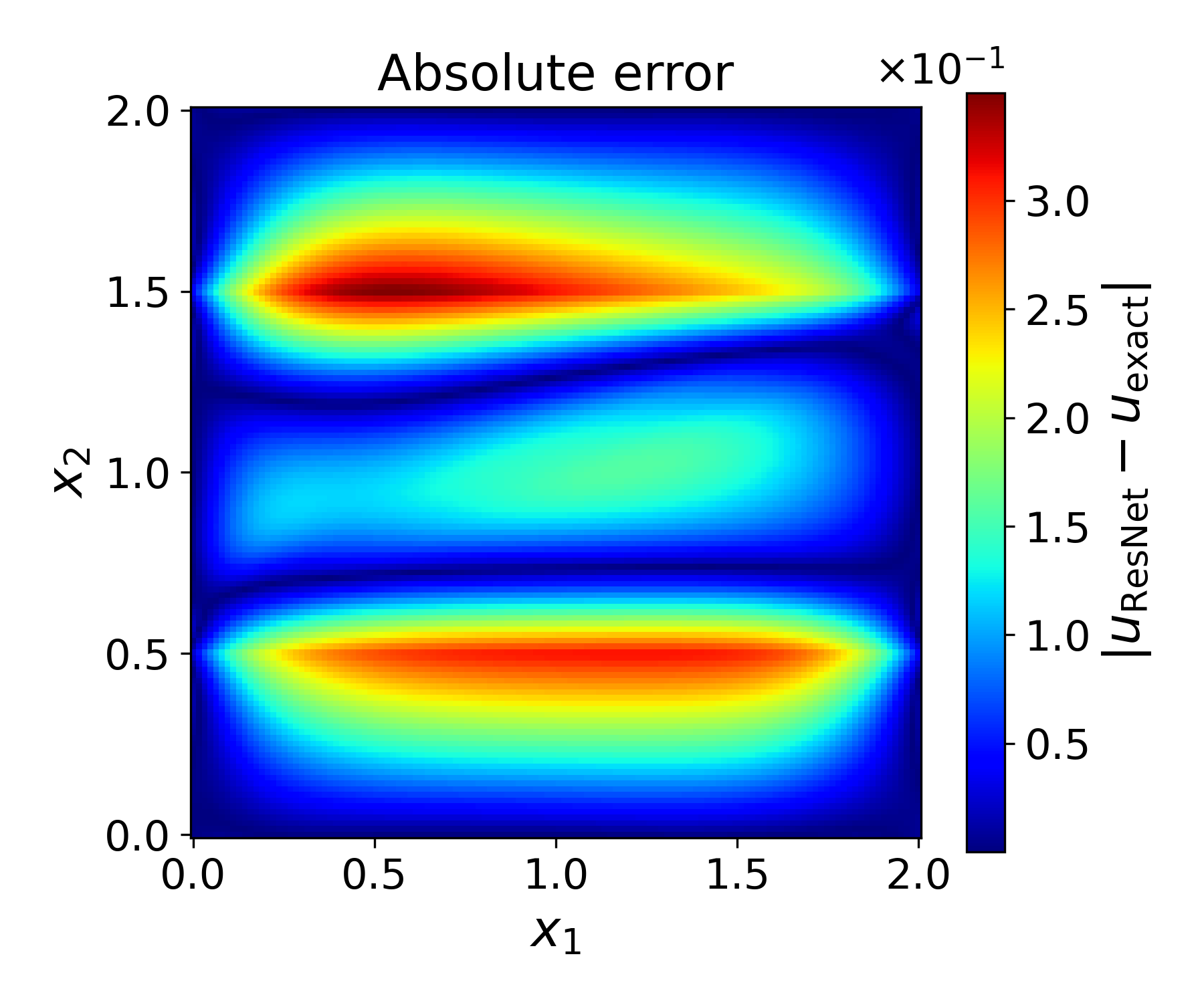}
    \caption{\(| u_{exact} - u_{\mathrm{ResNet}}|\)}
    \end{subfigure}
    \begin{subfigure}[t]{0.32\textwidth}
    \centering
    \includegraphics[width = \linewidth]{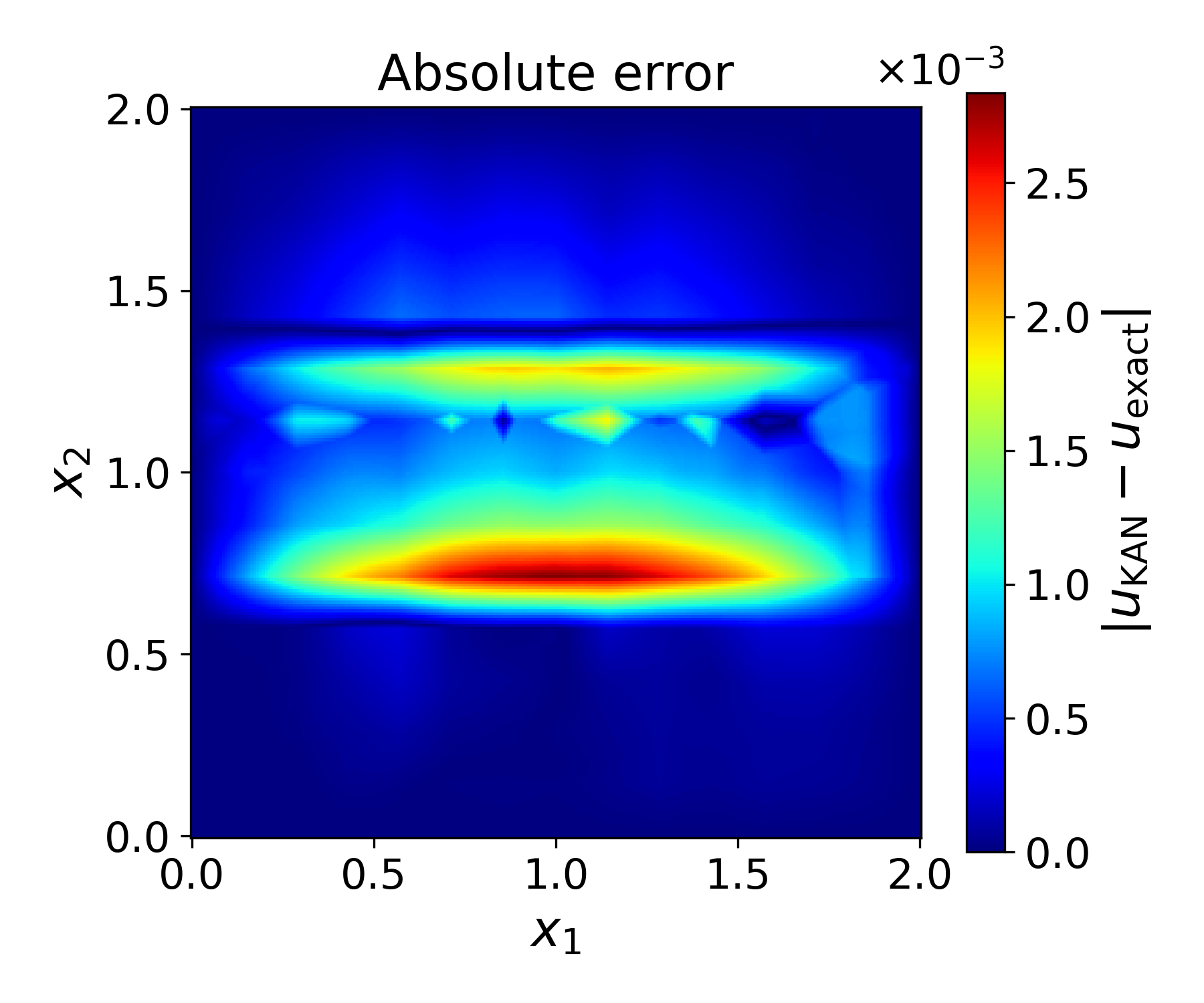}
    \caption{\(|u_{exact} - u_{\mathrm{KAN}}|\)}
    \end{subfigure}\hfill
    \vspace{0.2em}
    \begin{subfigure}[t]{0.45\textwidth}
    \centering
    \includegraphics[width = \linewidth]{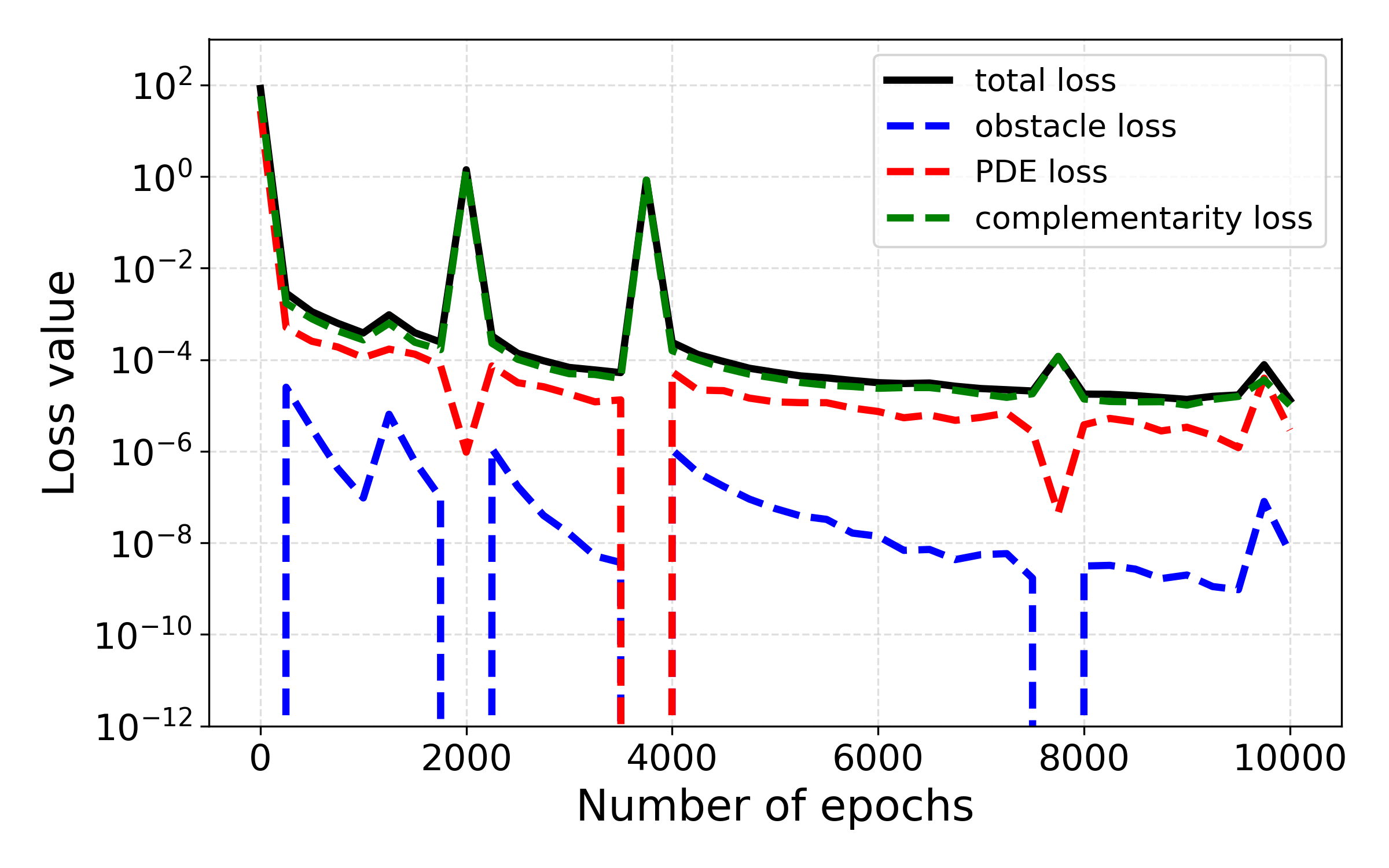}
    \caption{Loss components vs.\ epochs}
    \end{subfigure}
    \begin{subfigure}[t]{0.45\textwidth}
    \centering
    \includegraphics[width = \linewidth]{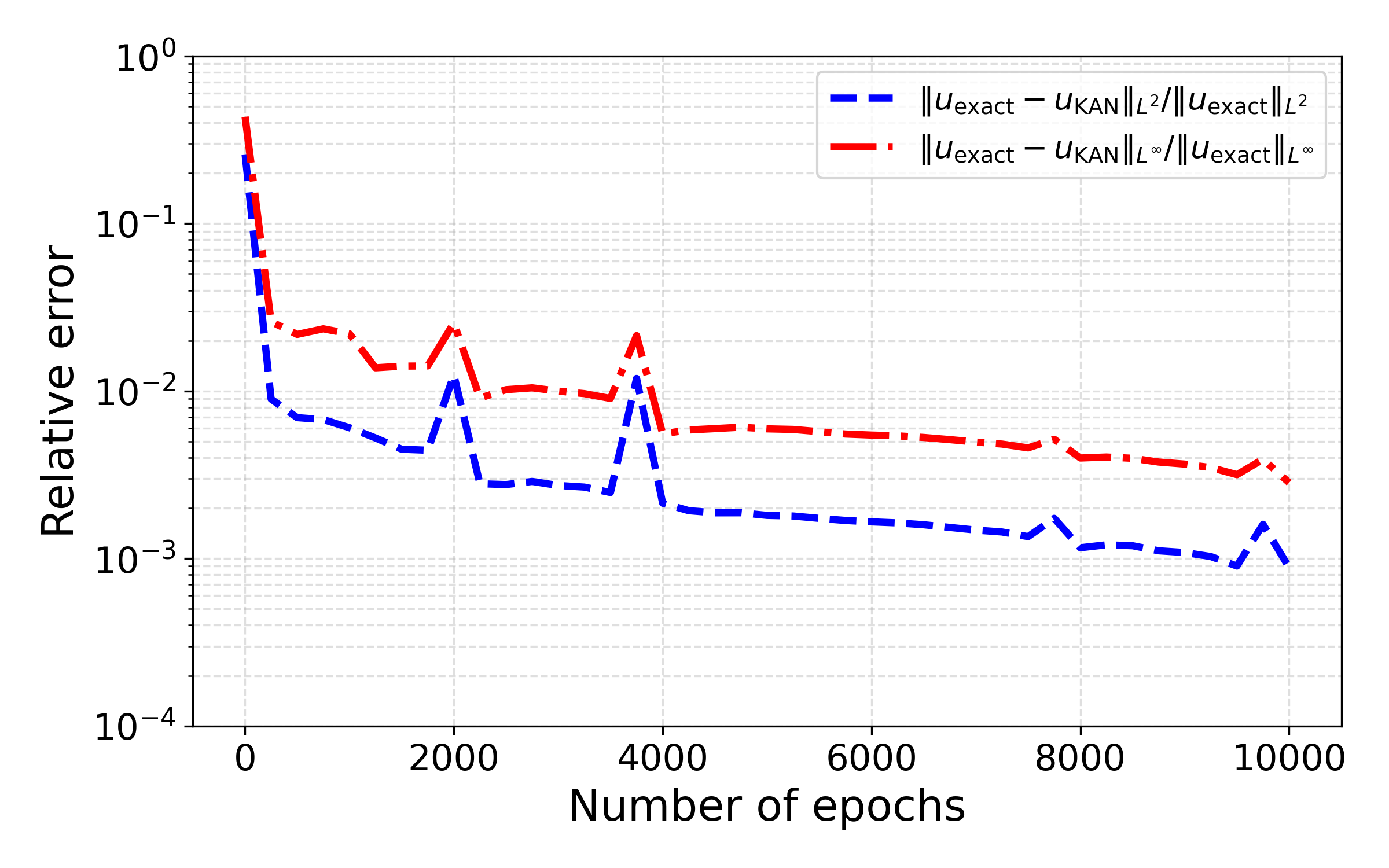}
    \caption{Relative error vs.\ epochs}
    \end{subfigure}\hfill
    \caption{\textbf{Comparison of PINN, ResNet, and KAN approximations for the two-dimensional \(p\)-Laplacian obstacle problem.}
Panels (a)--(c) show the obstacle function, exact solution, and \gls{kan} approximation. Panels (d)--(f) compare the absolute pointwise errors of the \gls{pinn}, \gls{resnet}, and \gls{kan} models. The largest \gls{kan} errors are concentrated near the free boundary. Panels (g) and (h) show the training loss components and relative errors over the training epochs. The final relative \(L^2\) and \(L^\infty\) errors of the \gls{kan} approximation are \(8.818\times10^{-4}\) and \(2.836\times10^{-3}\), respectively.}
\label{fig:plop_result}
\end{figure}

For the nonlinear \gls{plop}, we consider the computational domain
\[
    \overline{\Omega}=[0,2]\times[0,2],
\]
and set \(p=4\). The \gls{plop} is
\begin{equation}
\begin{cases}
    u(\vec{x}) \geq \psi(\vec{x}), & \vec{x}\in\Omega,\\
    -\nabla\cdot\left(|\nabla u(\vec{x})|^{p-2}\nabla u(\vec{x})\right)+1 \geq 0,
    & \vec{x}\in\Omega,\\
    \bigl(u(\vec{x})-\psi(\vec{x})\bigr)
    \left[
    -\nabla\cdot\left(|\nabla u(\vec{x})|^{p-2}\nabla u(\vec{x})\right)+1
    \right]=0,
    & \vec{x}\in\Omega,\\
    u(\vec{s})=g(\vec{s}), & \vec{s}\in\partial\Omega.
\end{cases}
\label{eq:appendix_plop_pde}
\end{equation}
Here, the nonlinear differential operator is the \(p\)-Laplacian,
\[
    -\Delta_p u
    =
    -\nabla\cdot\left(|\nabla u|^{p-2}\nabla u\right).
\]
We use the following piecewise obstacle function:
\begin{equation}
    \psi(x_1,x_2)
    =
    \begin{cases}
    1, & 0.5\leq x_1\leq 1.5,\\
    0, & \text{otherwise}.
    \end{cases}
    \label{eq:appendix_plop_obstacle}
\end{equation}
The exact solution used for evaluation is
\begin{equation}
    u_{\mathrm{exact}}(x_1,x_2)
    =
    \begin{cases}
    \dfrac{3}{4}\left|x_1+7.75086\right|^{4/3}-11.50434,
    & x_1<0.5,\\[0.9em]
    1,
    & 0.5\leq x_1\leq 1.5,\\[0.9em]
    \dfrac{3}{4}\left|-x_1+9.75086\right|^{4/3}-11.50434,
    & x_1>1.5.
    \end{cases}
    \label{eq:appendix_plop_exact_solution}
\end{equation}
The Dirichlet boundary data are prescribed from the exact solution:
\[
    g(\vec{s})=u_{\mathrm{exact}}(\vec{s}),
    \qquad \vec{s}\in\partial\Omega.
\]
Figure~\ref{fig:plop_result} presents the obstacle function, exact solution, \gls{kan} approximation, pointwise errors, and training history.

\FloatBarrier

Figure~\ref{fig:plop_result} shows that the \gls{kan} approximation closely agrees with the exact nonlinear \(p\)-Laplacian obstacle solution. The remaining error is concentrated mainly near the free boundary, where the solution transitions between the contact and non-contact regions. The final relative \(L^2\) and \(L^\infty\) errors are \(8.818\times10^{-4}\) and \(2.836\times10^{-3}\), respectively.

The ablation study in Figure~\ref{fig_params_gridsize_study}(b) shows a clear interaction between spline grid size and network width. For the width-8 model, the relative \(L^2\) error is approximately \(1.9\times10^{-5}\) at \(G=5\), \(1.0\times10^{-4}\) at \(G=10\), \(9.5\times10^{-6}\) at \(G=15\), and \(1.8\times10^{-4}\) at \(G=20\). For the width-48 model, the corresponding errors are approximately \(4.1\times10^{-4}\), \(8.4\times10^{-4}\), \(8.818\times10^{-4}\), and \(1.4\times10^{-3}\). The smallest error in this ablation is obtained with the width-8 model at \(G=15\). These results show that increasing spline resolution or network width does not necessarily improve accuracy.

Having considered the stationary obstacle problems, we next study the time-dependent \gls{sp} in Subsection~\ref{exp:time-dependent-stefan-problem}. In this case, the temperature field and free boundary evolve simultaneously in time.

\FloatBarrier
\subsection{Time-dependent Stefan Problem}\label{exp:time-dependent-stefan-problem}
For the time-dependent experiment, we consider a two-dimensional one-phase \gls{sp} based on \cite{beckett2001moving,bansch2013ale,wang2021deep}. The unknown temperature field is denoted by \(u(x_1,x_2,t)\), and the moving free boundary is represented by
\[
    x_1=s(x_2,t).
\]
We consider
\[
    0\leq x_1\leq 2.25,
    \qquad
    0\leq x_2\leq 1,
    \qquad
    0\leq t\leq 1.
\]
The moving spatial domain is
\begin{equation}
    \Omega(t)
    =
    \left\{
    (x_1,x_2)\in\Omega:
    0<x_1<s(x_2,t),\ 0<x_2<1
    \right\},
    \qquad t\in[0,1].
    \label{eq:appendix_stefan_domain}
\end{equation}
The one-phase \gls{sp} is to find \(u(x_1,x_2,t)\) and \(s(x_2,t)\) such that
\begin{equation}
\begin{cases}
    \partial_t u(x_1,x_2,t)-\Delta u(x_1,x_2,t)=0,
    & (x_1,x_2)\in\Omega(t),\ t\in(0,1],\\[0.3em]
    u(s(x_2,t),x_2,t)=0,
    & 0<x_2<1,\ t\in(0,1],\\[0.3em]
    \partial_t s(x_2,t)
    +
    \partial_{x_1}u(s(x_2,t),x_2,t)
    -
    \partial_{x_2}s(x_2,t)\,
    \partial_{x_2}u(s(x_2,t),x_2,t)
    =0,
    & 0<x_2<1,\ t\in(0,1],\\[0.3em]
    u(x_1,x_2,0)=u_0(x_1,x_2),
    & (x_1,x_2)\in\Omega(0),\\[0.3em]
    s(x_2,0)=s_0(x_2),
    & 0<x_2<1.
\end{cases}
\label{eq:appendix_stefan_problem}
\end{equation}
The first equation is the heat equation in the moving domain. The second and third equations impose the interface temperature and Stefan conditions, respectively. The final two equations prescribe the initial temperature field and initial interface location.

For this problem, the exact moving free boundary is prescribed by
\begin{equation}
    s_{\mathrm{exact}}(x_2,t)
    =
    \frac{1}{2}x_2+\frac{5}{4}t+\frac{1}{2},
    \label{eq:appendix_stefan_exact_interface}
\end{equation}
and the exact temperature solution is
\begin{equation}
    u_{\mathrm{exact}}(x_1,x_2,t)
    =
    \exp\left(
    \frac{5}{4}t-x_1+\frac{1}{2}x_2+\frac{1}{2}
    \right)-1.
    \label{eq:appendix_stefan_exact_solution}
\end{equation}
Therefore, the initial conditions are
\begin{equation}
    u_0(x_1,x_2)
    =
    \exp\left(
    -x_1+\frac{1}{2}x_2+\frac{1}{2}
    \right)-1,
    \qquad
    s_0(x_2)=\frac{1}{2}x_2+\frac{1}{2}.
    \label{eq:appendix_stefan_initial_conditions}
\end{equation}
Following Algorithm~\ref{alg:stefan_solver}, one \gls{kan} approximates the temperature field \(u_{\theta}(x_1,x_2,t)\), and a second \gls{kan} approximates the moving interface \(\hat{s}_{\eta}(x_2,t)\). The temperature network takes \((x_1,x_2,t)\) as input, whereas the interface network takes \((x_2,t)\) as input. The fixed-boundary data on \(x_1=0\), \(x_2=0\), and \(x_2=1\) are imposed through the construction of the temperature model. The initial conditions for both \(u\) and \(s\) are enforced through \(\mathcal{L}_{\mathrm{ic}}\). The remaining loss terms enforce the heat equation, interface temperature condition, and Stefan free-boundary condition:
\begin{equation}
    \mathcal{L}_{\mathrm{SP}}
    =
    \lambda_{\mathrm{heat}}\mathcal{L}_{\mathrm{heat}}
    +
    \lambda_{\mathrm{ic}}\mathcal{L}_{\mathrm{ic}}
    +
    \lambda_{\mathrm{int}}\mathcal{L}_{\mathrm{int}}
    +
    \lambda_{\mathrm{fb}}\mathcal{L}_{\mathrm{fb}}.
    \label{eq:appendix_stefan_total_loss}
\end{equation}
Here, \(\mathcal{L}_{\mathrm{heat}}\) enforces the heat equation, \(\mathcal{L}_{\mathrm{ic}}\) penalizes the initial conditions, \(\mathcal{L}_{\mathrm{int}}\) enforce the interface temperature condition, and \(\mathcal{L}_{\mathrm{fb}}\) enforce the Stefan free-boundary condition.
\FloatBarrier
\begin{figure}[p]
    \centering
    \setlength{\tabcolsep}{1pt}
    \renewcommand{\arraystretch}{0.55}

    \begin{tabular}{ccc}

    \includegraphics[width=0.295\textwidth]{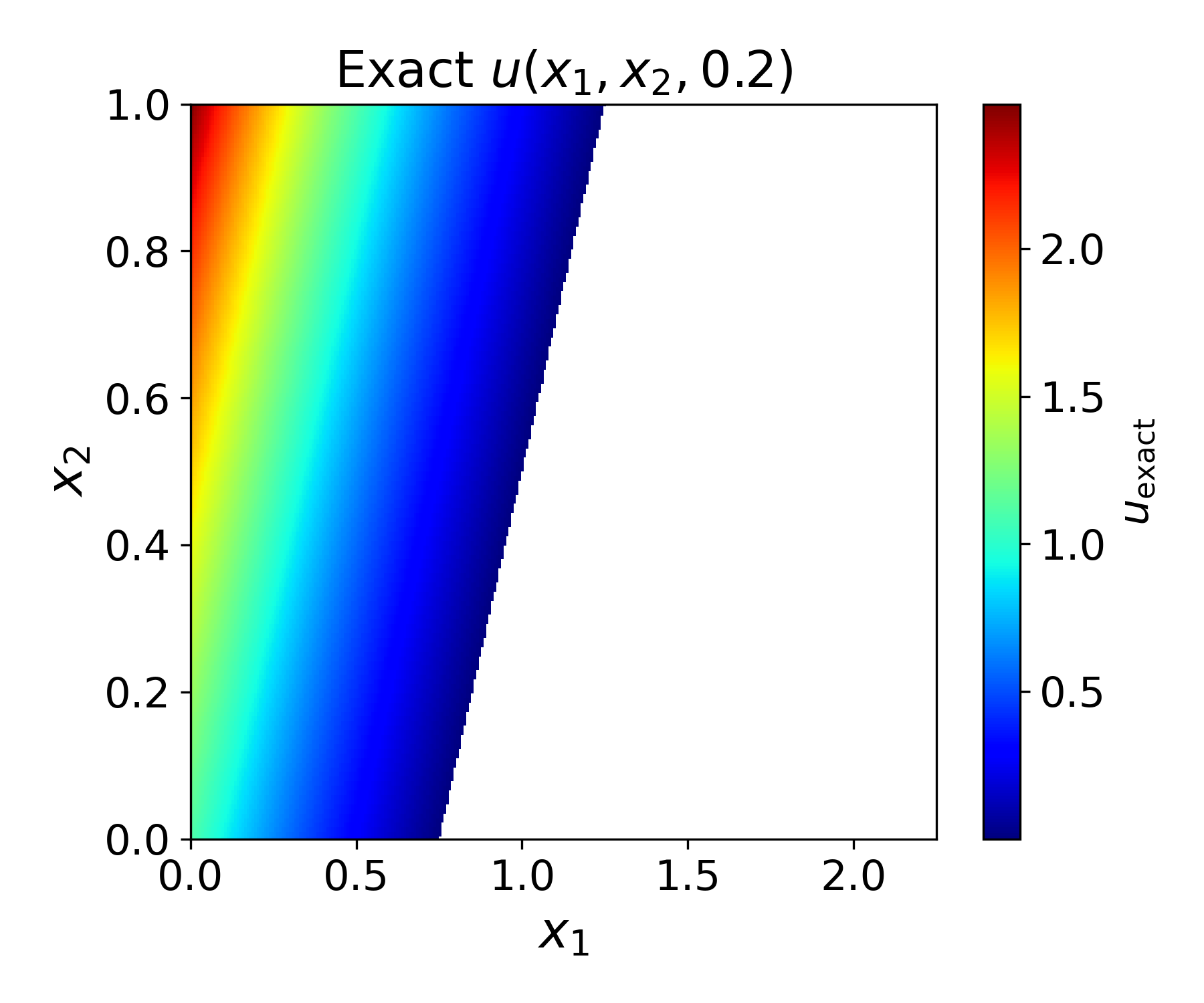}
    &
    \includegraphics[width=0.295\textwidth]{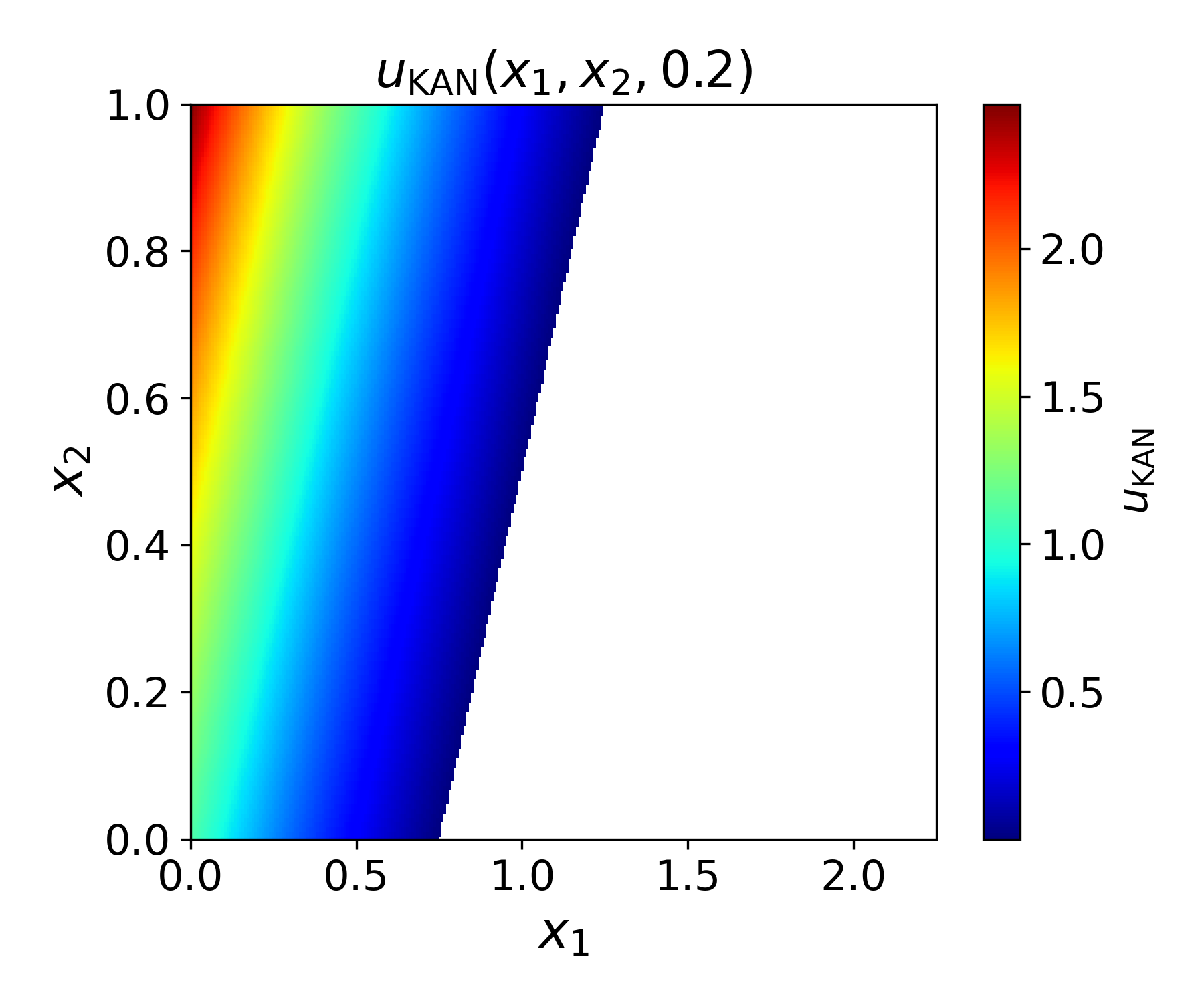}
    &
    \includegraphics[width=0.295\textwidth]{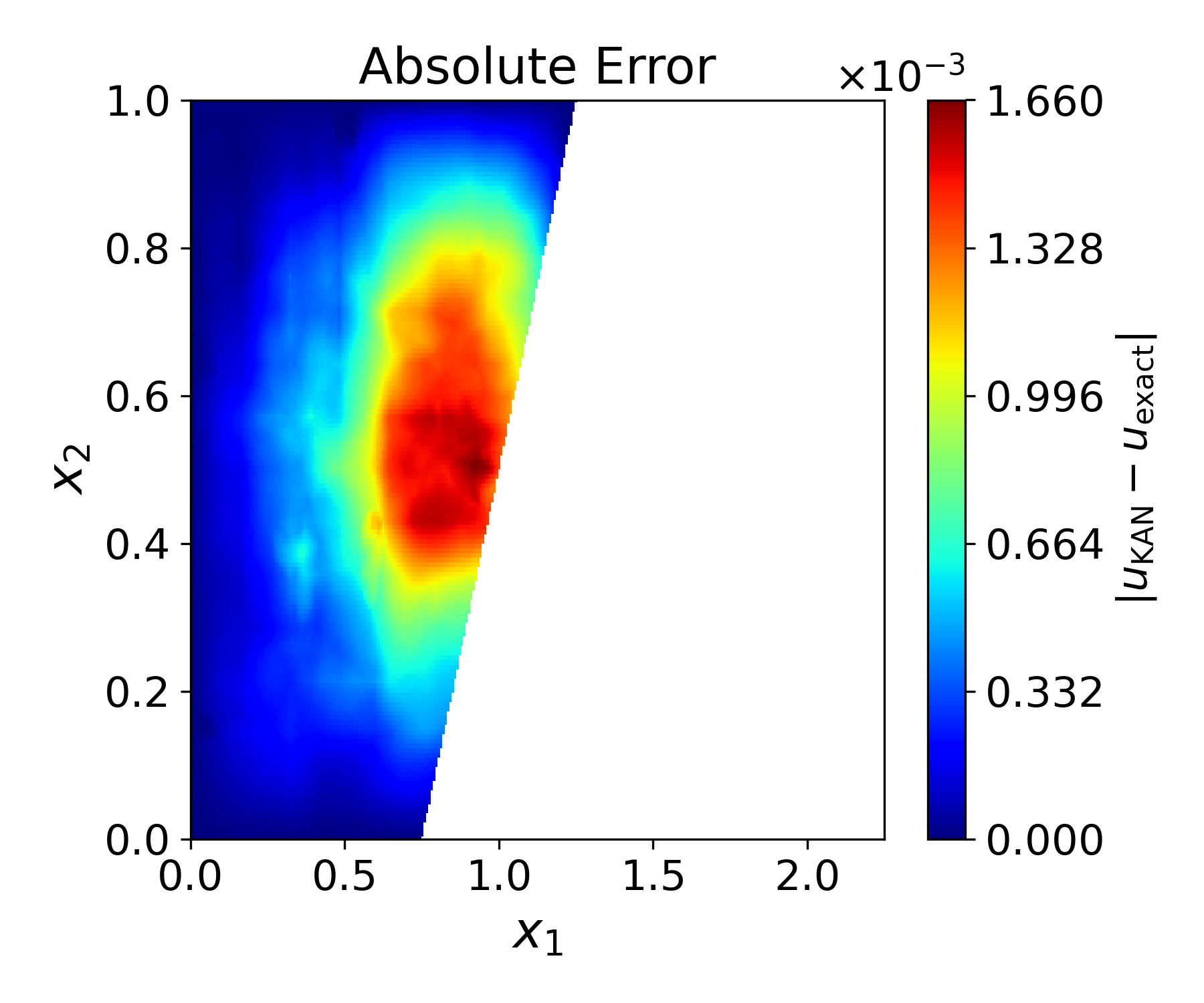}
    \\[-0.4em]

    \includegraphics[width=0.295\textwidth]{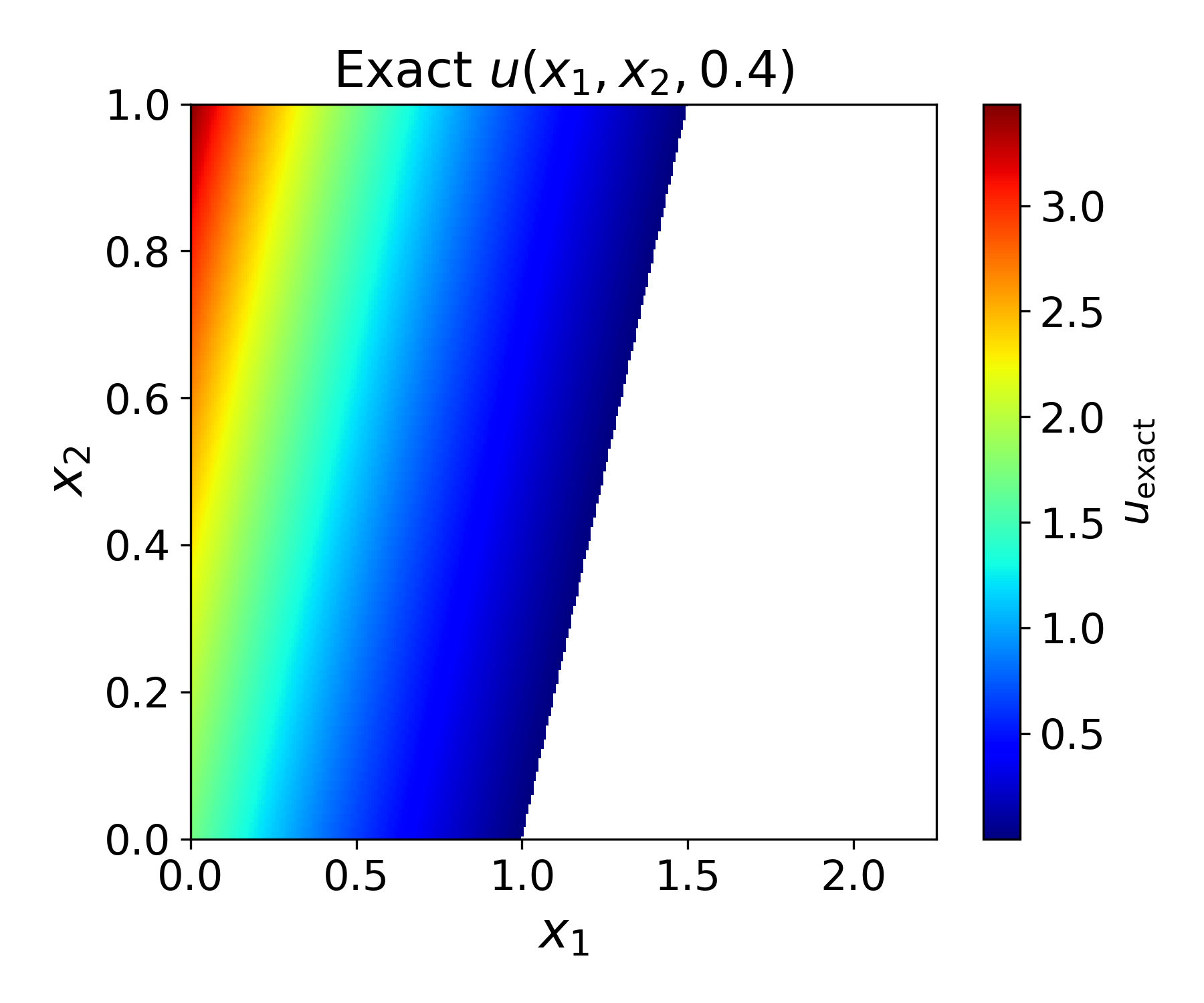}
    &
    \includegraphics[width=0.295\textwidth]{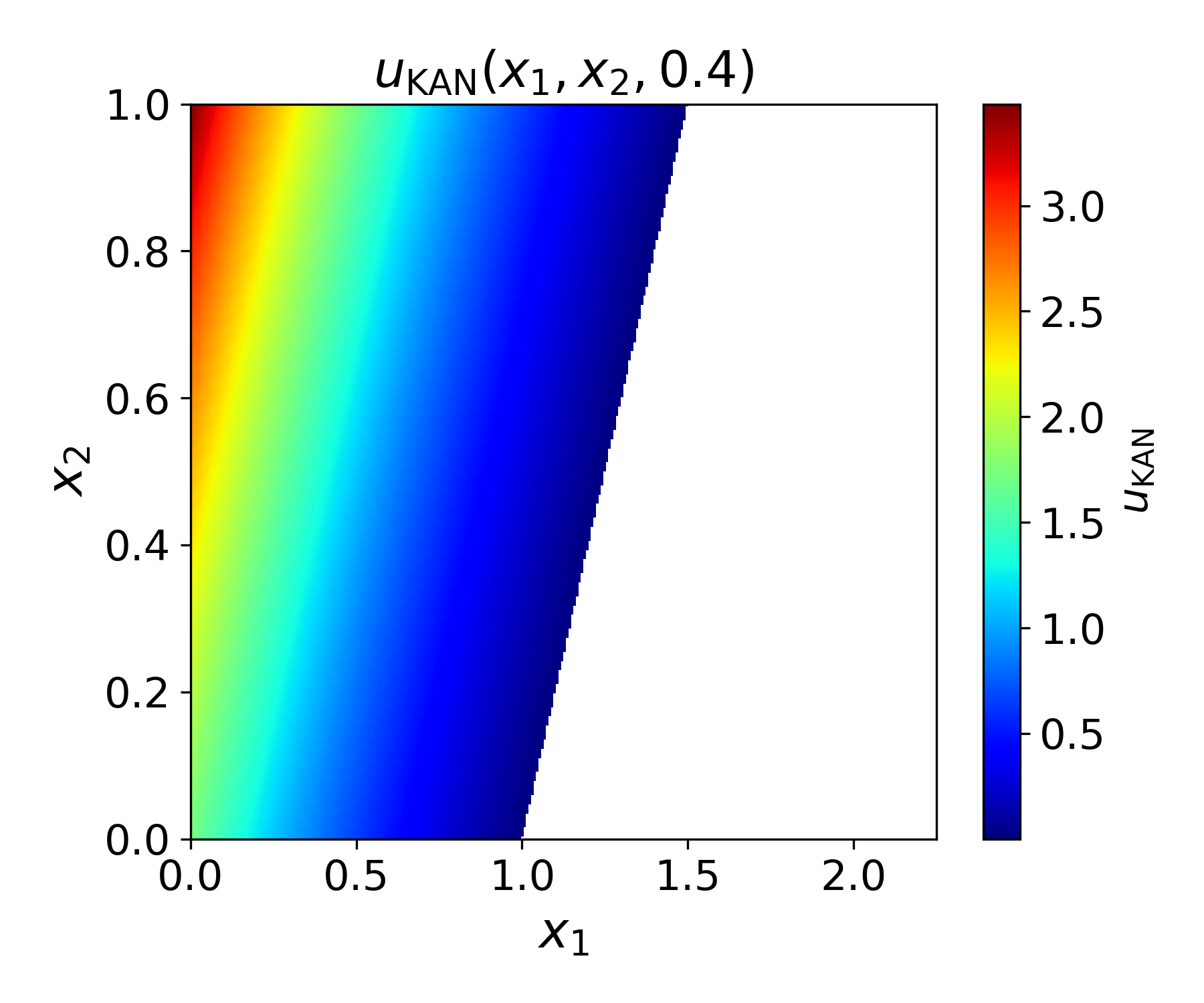}
    &
    \includegraphics[width=0.295\textwidth]{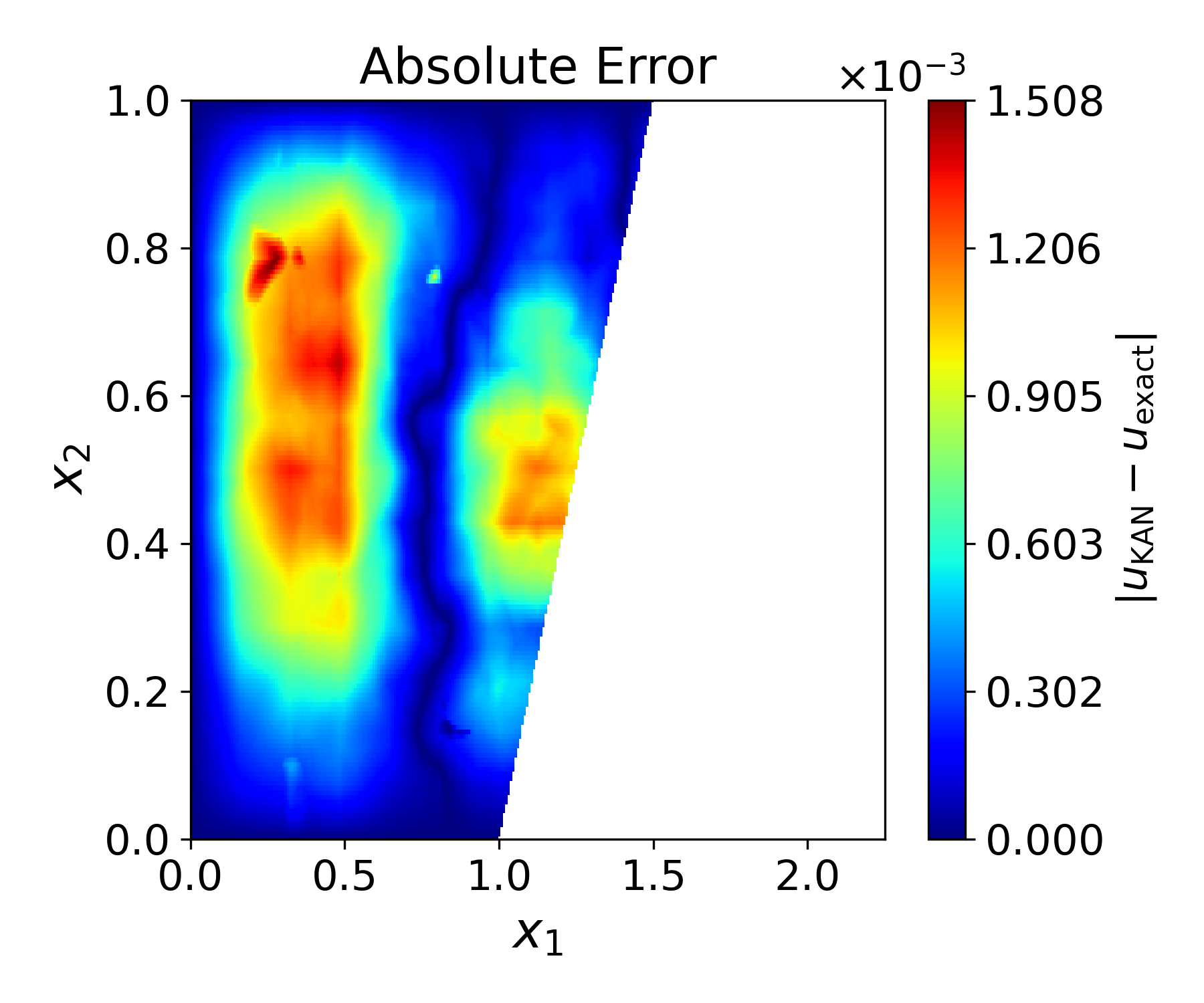}
    \\[-0.4em]

    \includegraphics[width=0.295\textwidth]{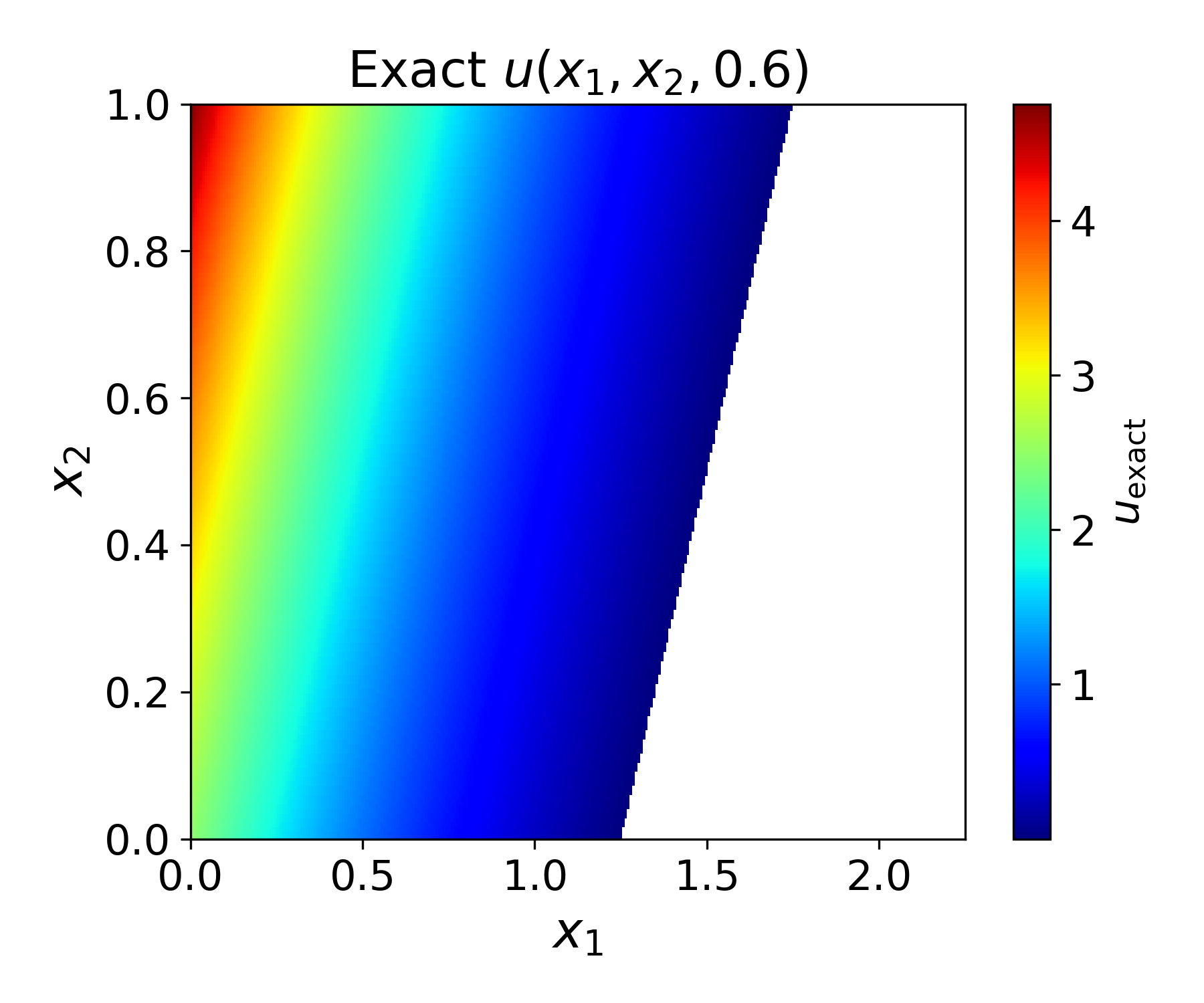}
    &
    \includegraphics[width=0.295\textwidth]{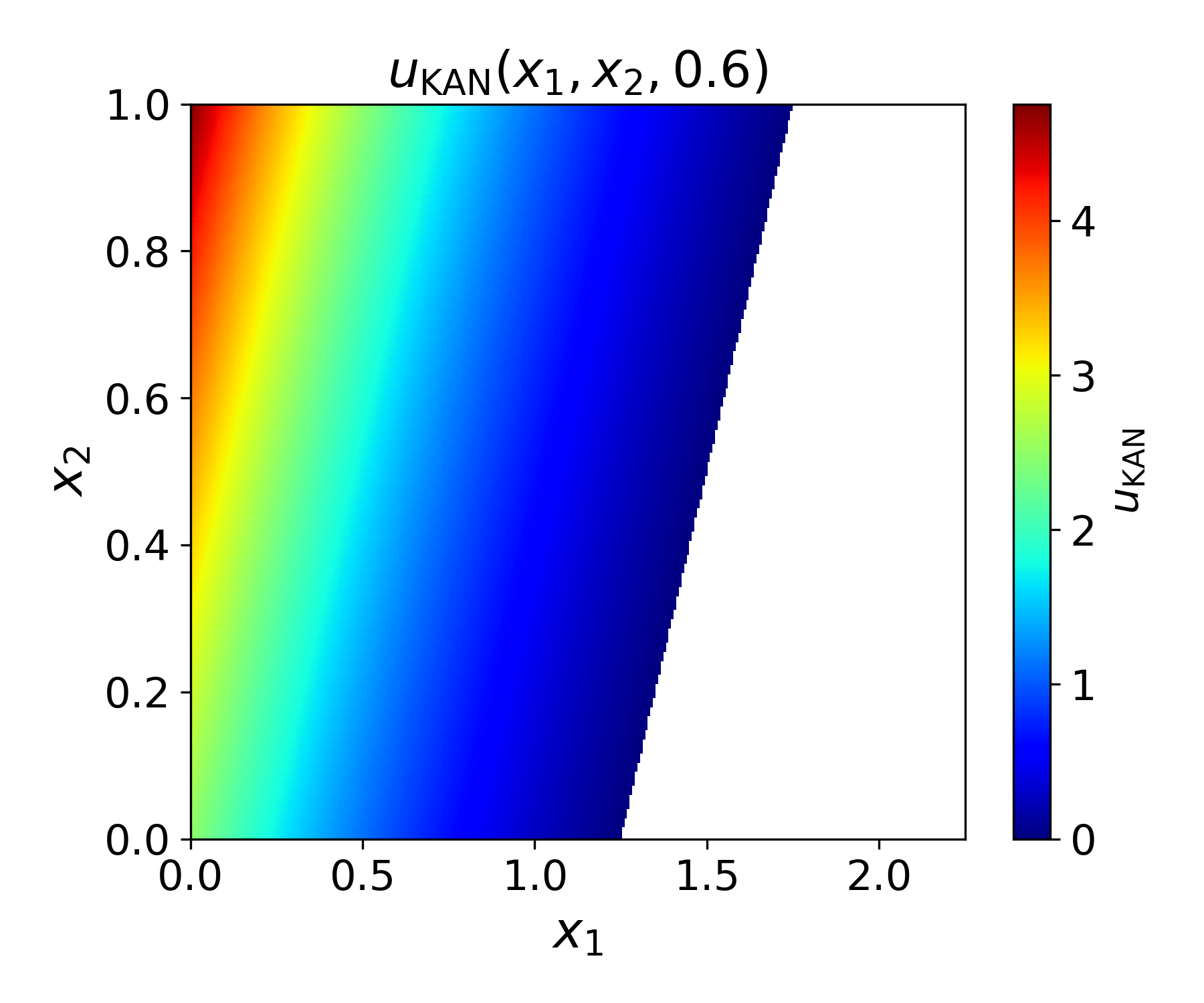}
    &
    \includegraphics[width=0.295\textwidth]{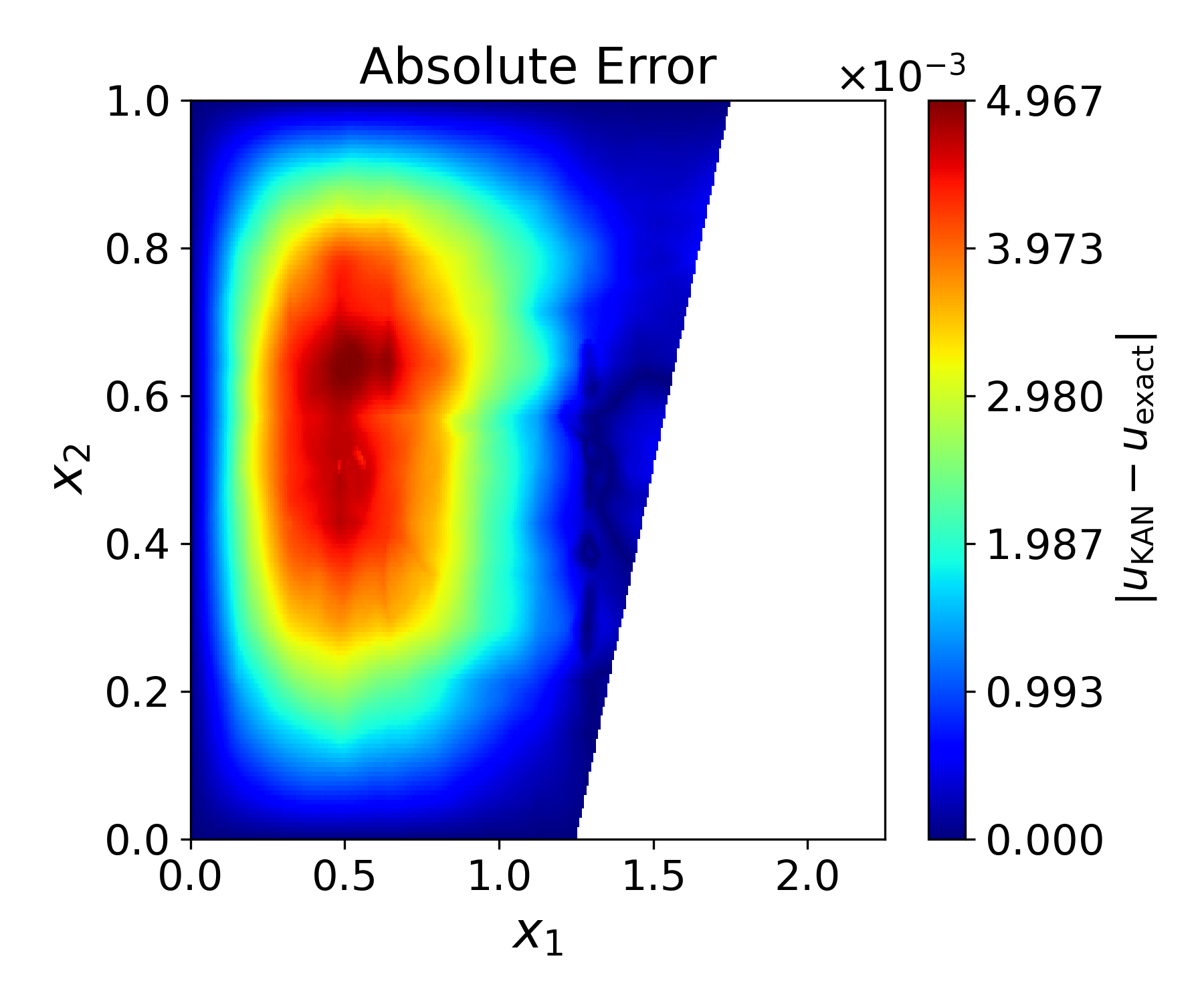}
    \\[-0.4em]

    \includegraphics[width=0.295\textwidth]{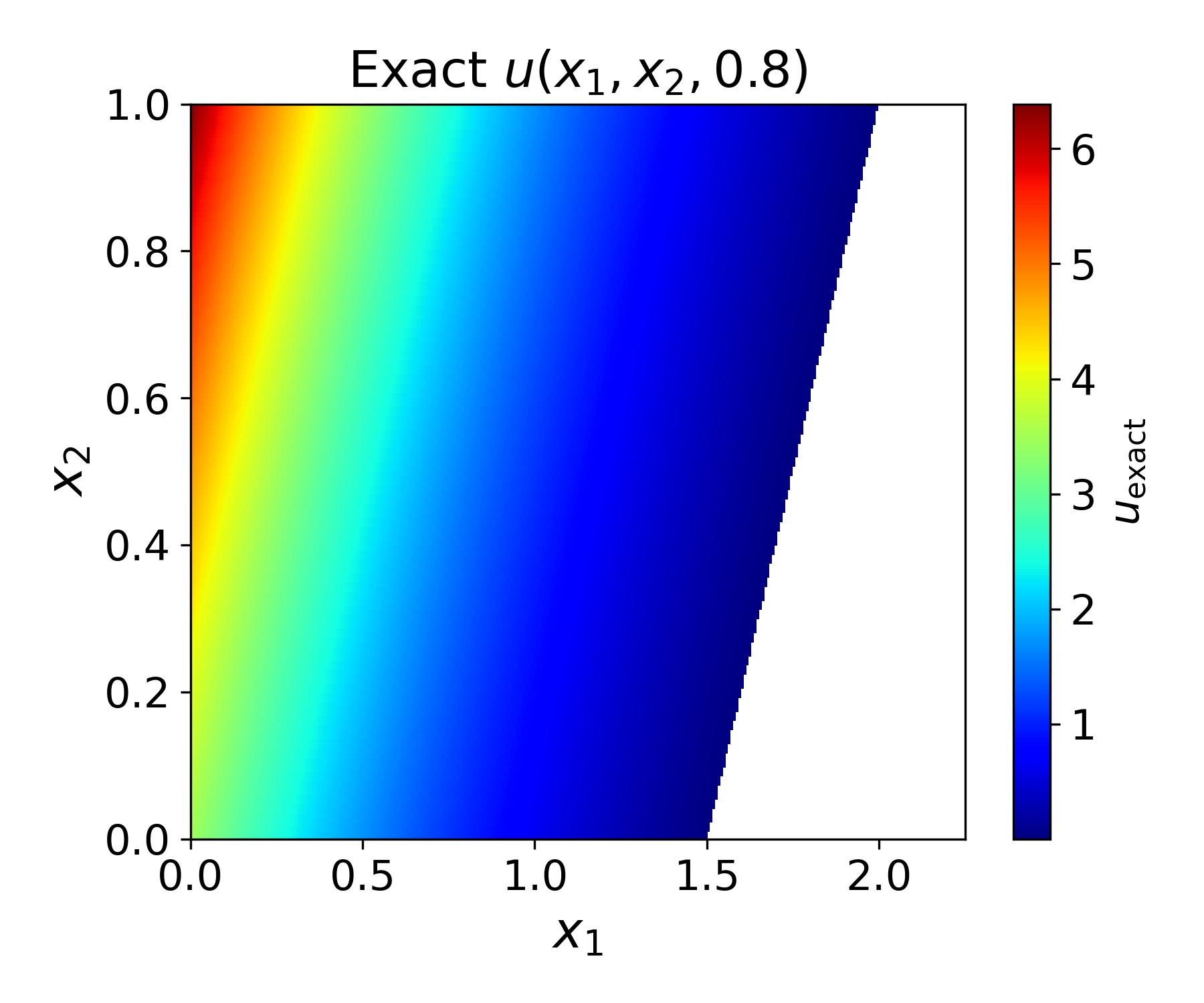}
    &
    \includegraphics[width=0.295\textwidth]{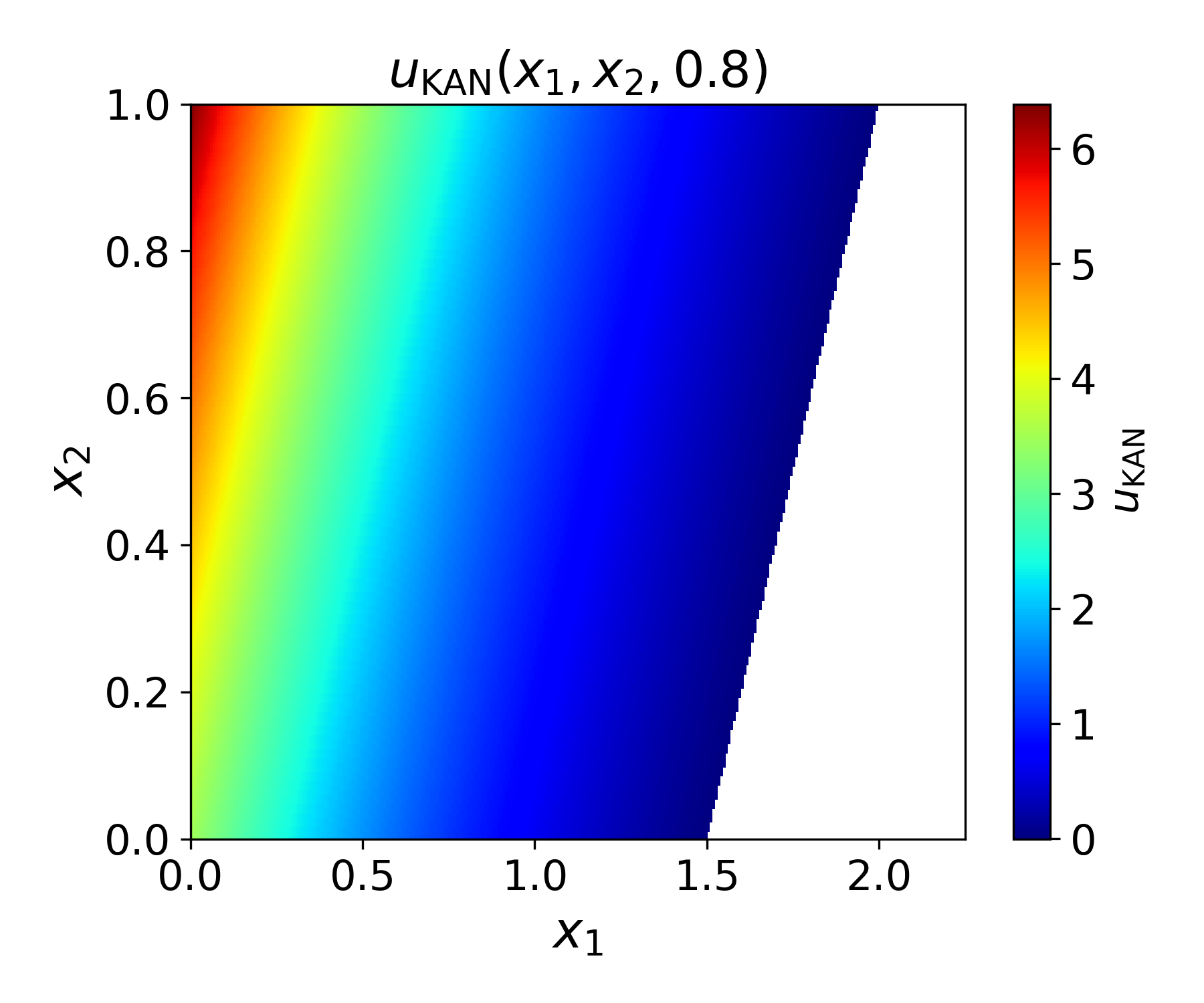}
    &
    \includegraphics[width=0.295\textwidth]{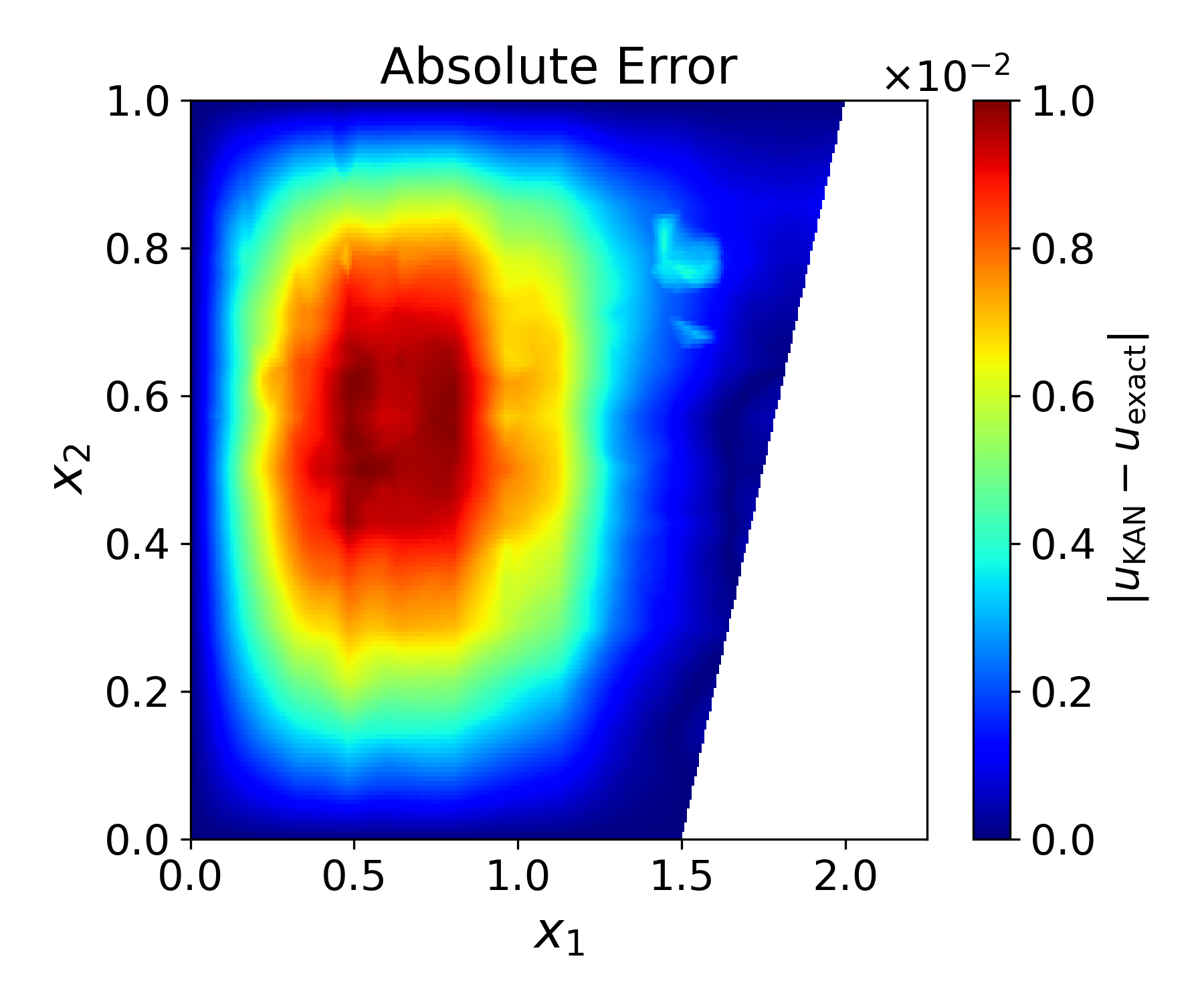}

    \end{tabular}

    \vspace{-0.8em}

    \caption{
    \textbf{KAN approximation for the two-dimensional one-phase Stefan problem at multiple time snapshots.}
    The four rows correspond to \(t=0.2,0.4,0.6,\) and \(0.8\). The columns show the exact solution
    \(u_{\mathrm{exact}}(x_1,x_2,t)\), the KAN approximation \(u_{\mathrm{KAN}}(x_1,x_2,t)\), and the absolute pointwise error
    \(|u_{\mathrm{KAN}}-u_{\mathrm{exact}}|\). The \gls{kan} model accurately approximates the temperature field on the moving domain, with the largest errors concentrated near the evolving interface.
    }
    \label{fig_stefan_problem_result}
\end{figure}
\FloatBarrier
\begin{figure}[p]
    \centering
    \begin{subfigure}[t]{0.32\textwidth}
    \centering
    \includegraphics[width = \linewidth]{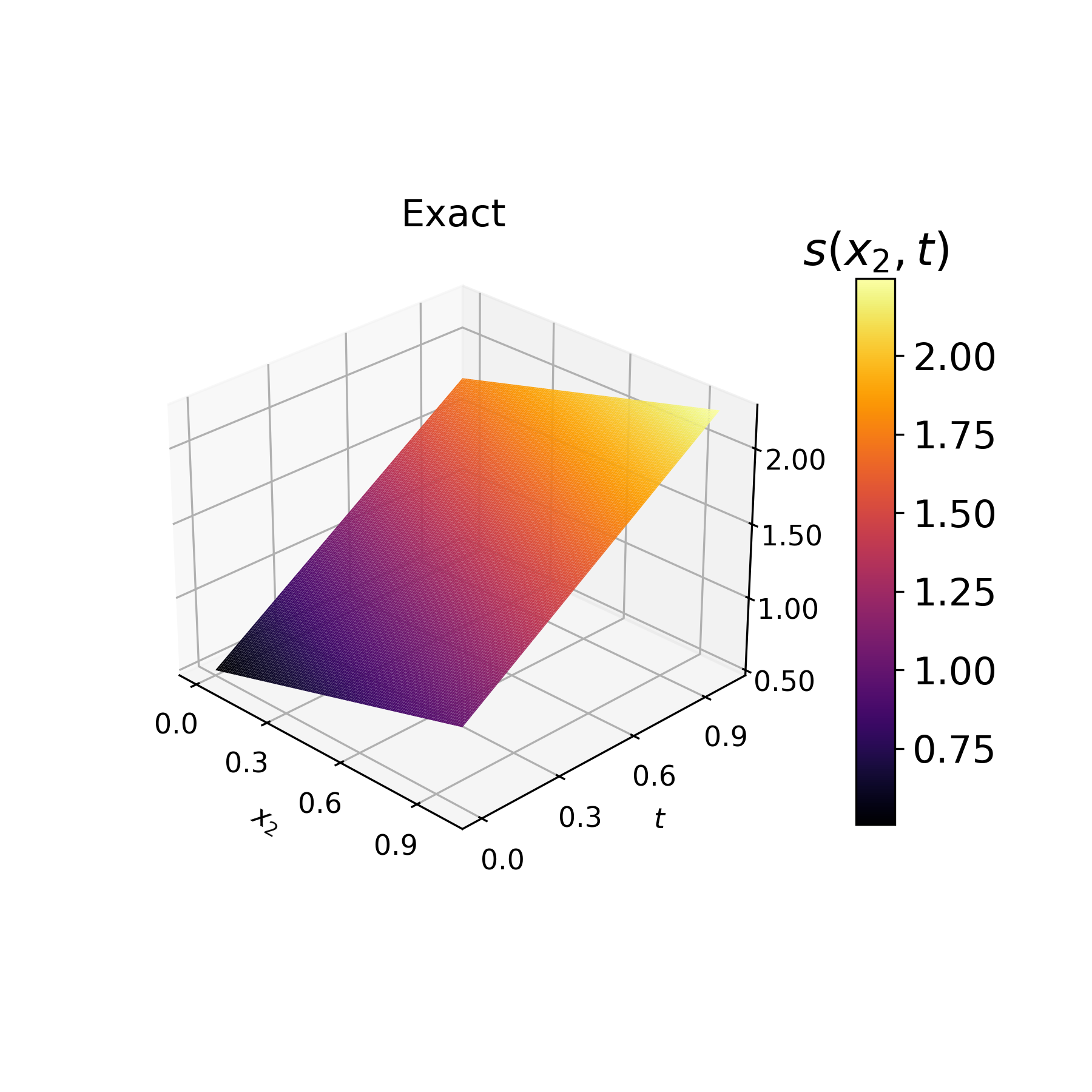}
    \caption{Exact moving free boundary}
    \end{subfigure}
    \begin{subfigure}[t]{0.32\textwidth}
    \centering
    \includegraphics[width = \linewidth]{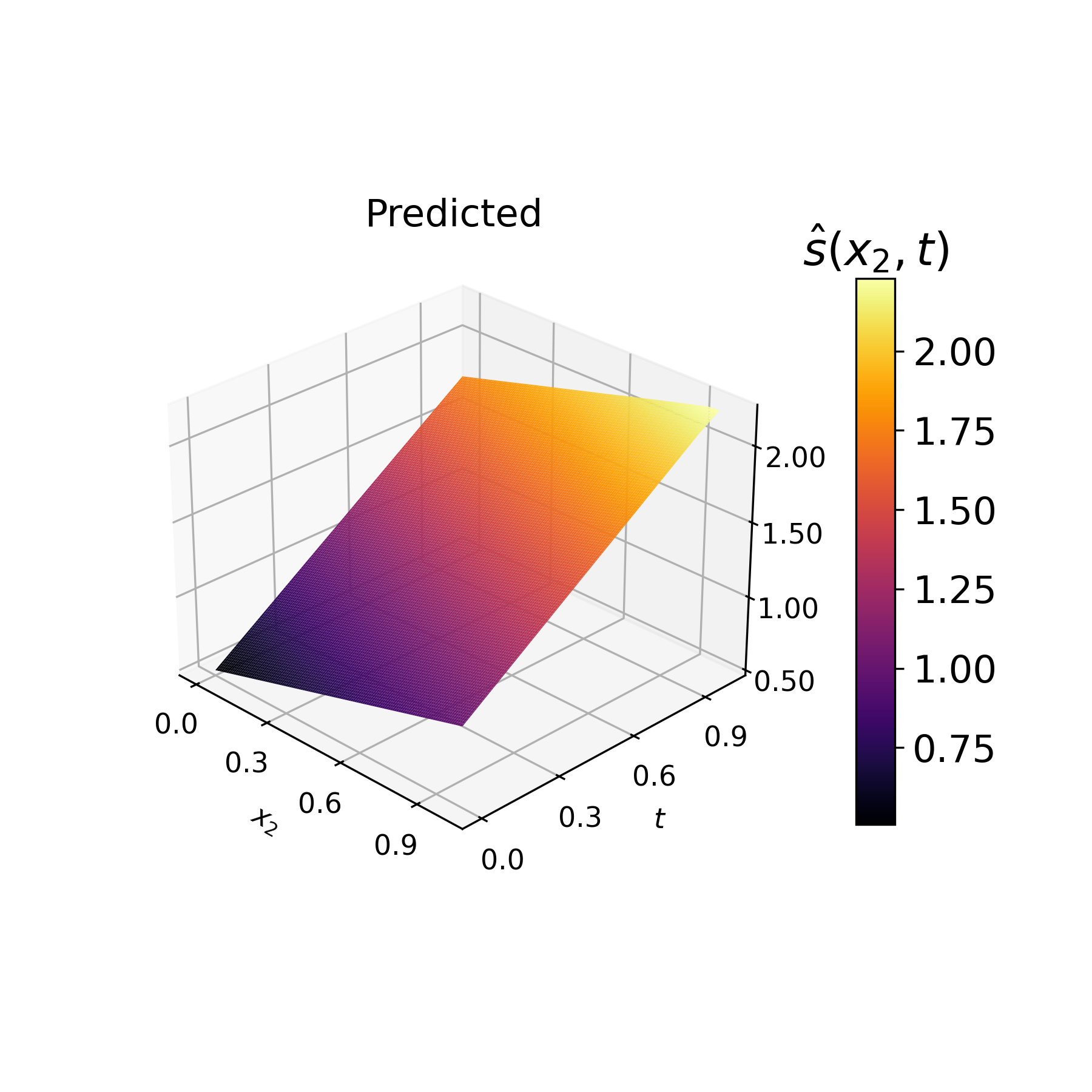}
    \caption{Predicted moving free boundary}
    \end{subfigure}
    \begin{subfigure}[t]{0.32\textwidth}
    \centering
    \includegraphics[width = \linewidth]{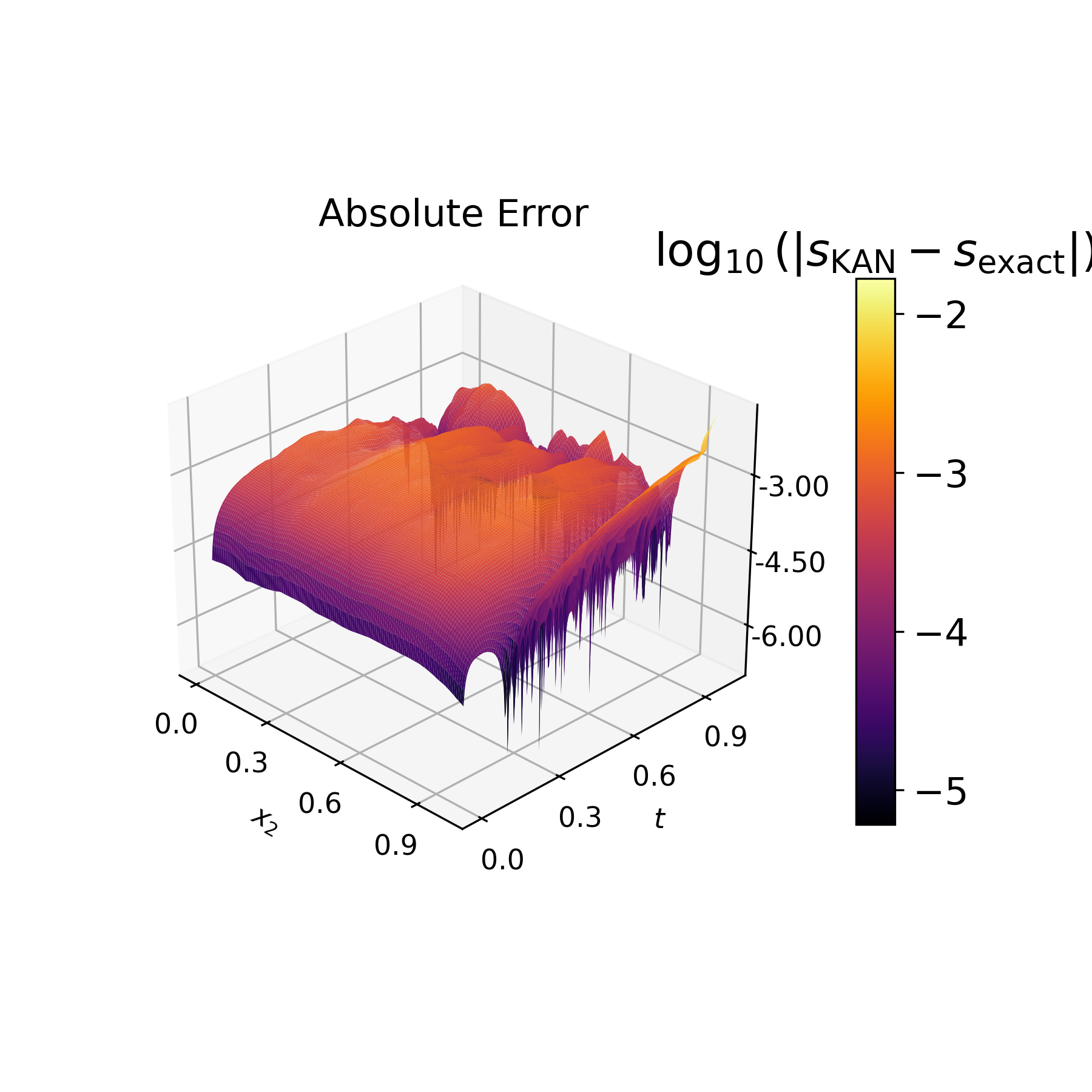}
    \caption{Logarithmic absolute error}
    \end{subfigure}\hfill
    \vspace{0.2em}
    \begin{subfigure}[t]{0.45\textwidth}
    \centering
    \includegraphics[width = \linewidth]{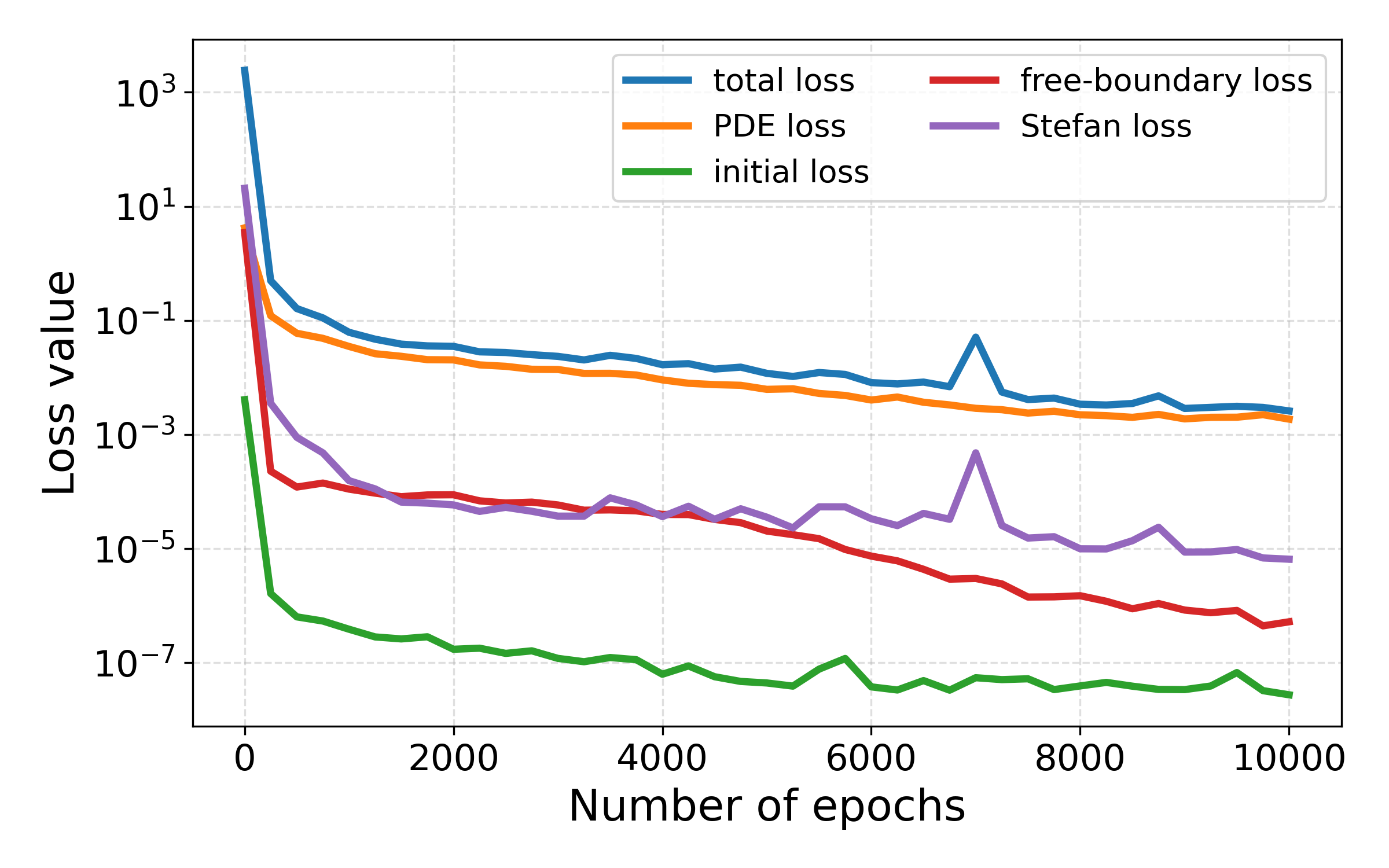}
    \caption{Loss components vs.\ epochs}
    \end{subfigure}
    \begin{subfigure}[t]{0.45\textwidth}
    \centering
    \includegraphics[width = \linewidth]{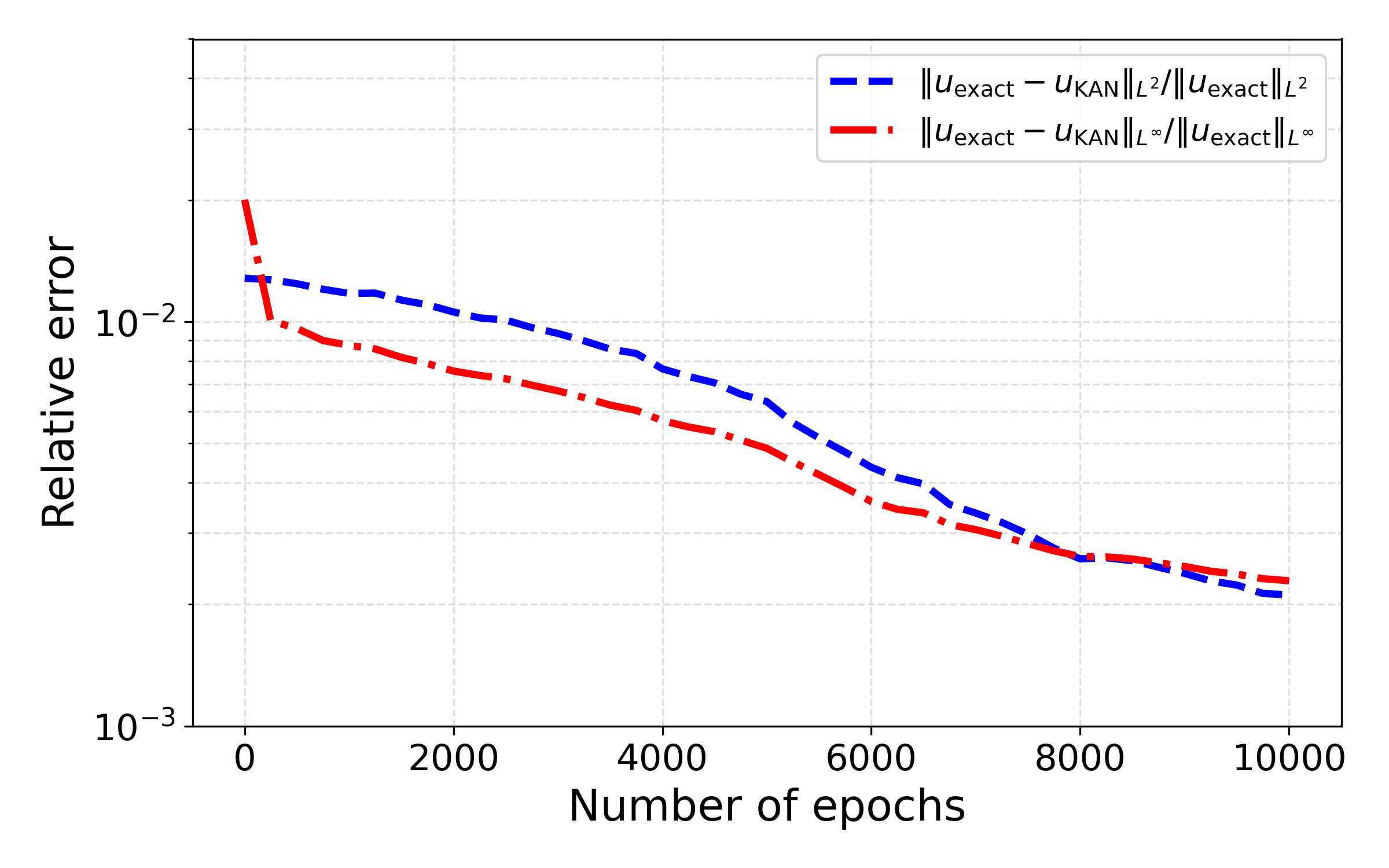}
    \caption{Relative error vs.\ epochs}
    \end{subfigure}
    \caption{\textbf{Moving free-boundary approximation and training history for the two-dimensional one-phase Stefan problem.} Panels (a) and (b) show the exact and predicted moving interfaces, respectively. Panel (c) shows the logarithmic absolute interface error, \(\log_{10}|s_{\mathrm{KAN}}-s_{\mathrm{exact}}|\). Panels (d) and (e) show the training loss components and relative solution errors over the training epochs. The final relative \(L^2\) error of the temperature field is \(5.966\times10^{-4}\).}
    \label{fig_stefan_problem_result_free_surface}
\end{figure}
\FloatBarrier
\begin{figure}[!htbp]
    \centering
    \includegraphics[width =\linewidth]{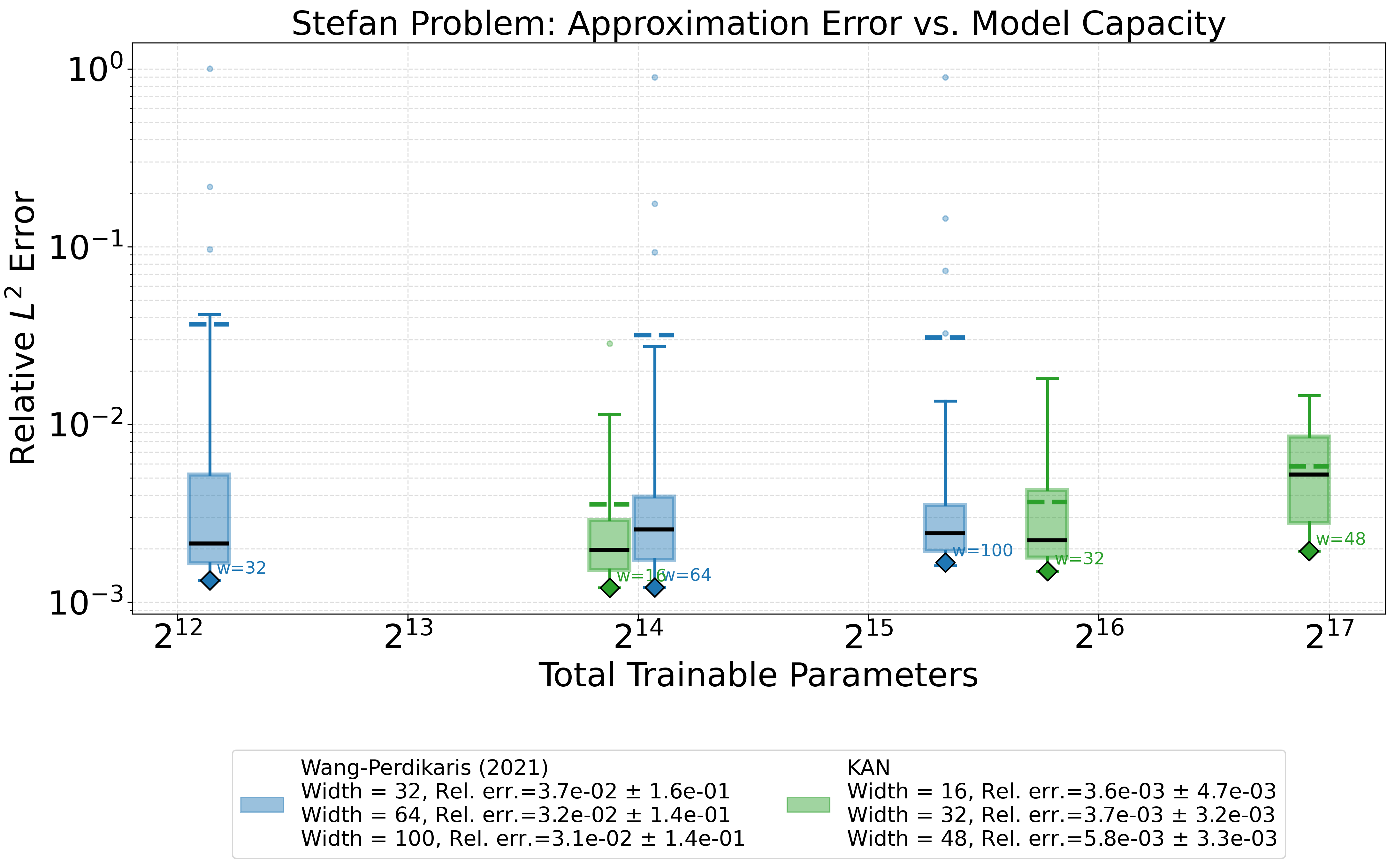}
    \caption{\textbf{Model-capacity comparison for the two-dimensional Stefan problem.}
Relative \(L^{2}\) error of the temperature field \(u\) versus the total number of trainable parameters for the baseline method in \cite{wang2021deep} and the proposed \gls{kan} models. The \gls{kan} models achieve lower errors with fewer trainable parameters, indicating improved parameter efficiency for this experiment.}
\label{fig:Stefan_params_study}
\end{figure}

Figures~\ref{fig_stefan_problem_result}--\ref{fig:Stefan_params_study} summarize the Stefan results. Figure~\ref{fig_stefan_problem_result} shows close agreement between \(u_{\mathrm{KAN}}\) and \(u_{\mathrm{exact}}\) at \(t=0.2,0.4,0.6,\) and \(0.8\). The pointwise errors are concentrated mainly near the moving interface. The maximum displayed error is of order \(10^{-3}\) at the earlier snapshots and reaches approximately an order of \(10^{-2}\) at \(t=0.8\).

Figure~\ref{fig_stefan_problem_result_free_surface} shows that the predicted interface \(\hat{s}(x_2,t)\) closely follows \(s_{\mathrm{exact}}(x_2,t)\) over the space-time domain. Panel (c) shows logarithmic interface errors mainly between \(-5\) and \(-2\), corresponding to absolute errors of approximately \(10^{-5}\) to \(10^{-2}\). The heat-equation, initial-condition, interface, and Stefan-condition losses all decrease during training. The final relative \(L^2\) error of the temperature field is \(5.966\times10^{-4}\).

Figure~\ref{fig:Stefan_params_study} compares model capacity. For the baseline method in \cite{wang2021deep}, widths \(32\), \(64\), and \(100\) give mean relative \(L^2\) errors of \(3.7\times10^{-2}\), \(3.2\times10^{-2}\), and \(3.1\times10^{-2}\), respectively. The corresponding standard deviations are \(1.6\times10^{-1}\), \(1.4\times10^{-1}\), and \(1.4\times10^{-1}\). The \gls{kan} models with widths \(16\), \(32\), and \(48\) give mean relative \(L^2\) errors of \(3.6\times10^{-3}\), \(3.7\times10^{-3}\), and \(5.8\times10^{-3}\), with standard deviations \(4.7\times10^{-3}\), \(3.2\times10^{-3}\), and \(3.3\times10^{-3}\), respectively. Thus, the tested \gls{kan} models reduce the mean error by about one order of magnitude relative to the baseline.

The spline-grid ablation in Figure~\ref{fig_params_gridsize_study}(c) shows the same non-monotonic behavior. At \(G=5\), the relative \(L^2\) errors for widths \(8\), \(16\), \(32\), and \(48\) are approximately \(6.1\times10^{-4}\), \(6.9\times10^{-4}\), \(5.1\times10^{-4}\), and \(5.8\times10^{-4}\), respectively. At \(G=20\), the corresponding errors increase to approximately \(5.8\times10^{-3}\), \(1.9\times10^{-3}\), \(8.0\times10^{-3}\), and \(1.1\times10^{-2}\). The width-32 model at \(G=5\) gives the smallest error among the configurations shown. These results show that network width and spline grid size should be selected jointly.
\FloatBarrier
\section{Conclusion}
\label{sec:conclusion}
\hspace{5mm}We developed a \gls{kan}-based physics-informed framework for free-boundary \glspl{pde}. The method represents the unknown solution using learnable univariate edge functions and enforces the governing equations and free-boundary constraints through residual-based losses. We evaluated the framework on the \gls{eop}, nonlinear \gls{plop}, and time-dependent one-phase \gls{sp}. The numerical results show that the proposed \gls{kan} models accurately approximate the solutions and resolve both stationary contact regions and moving interfaces. In the tested examples, the \gls{kan} models achieve lower relative errors and improved parameter efficiency compared with the \gls{pinn} and \gls{resnet} baselines. The ablation studies also show that increasing the spline grid size or network width does not necessarily improve accuracy. Future work will consider adaptive spline refinement, more complex geometries, higher-dimensional free-boundary problems, and theoretical error analysis.


\paragraph*{Data and code availability.}
No external datasets were used in this study. All data analyzed in the numerical experiments were generated from the computational procedures described in this paper. The source code and scripts used for data generation, model training, and analysis are publicly available in the GitHub repository at \url{https://github.com/phuongdongle/KAN_obstacle}.

\paragraph*{Research funding.}
This research received no specific grant from any funding agency in the public, commercial, or not-for-profit sectors.

\paragraph*{Conflict of interest.}
The author declares no conflict of interest.

\pagebreak
\bibliographystyle{plainnat}
\bibliography{references_paper}
\pagebreak


\appendix

\section{Acronyms and General Notations}
\label{appendix:abbreviation_notations}
This appendix summarizes the abbreviations and mathematical notation used throughout the paper. Table~\ref{tab:general_notations} lists the principal symbols used in the problem formulations, loss functions, and numerical experiments.
\printglossary[type=\acronymtype,toctitle={List of Abbreviations}]
\subsection*{List of General Notations}
\begin{table}[!htbp]
\centering
\caption{Notations in the paper.}
\label{tab:general_notations}
\begin{tabular}{c  c }
\toprule
\textbf{Notation} & \textbf{Description} \\
\midrule
\(\Omega\) & Spatial computational domain \\
\(\partial\Omega\) & Boundary of the spatial domain \\
\(\vec{x}\in\Omega\) & Spatial variable \\
\(t\in[0,T]\) & Time variable \\
\(u_{\mathrm{exact}}\) & Exact solution \\
\(u_{\theta}\) & \gls{kan} approximation of the solution \\
\(\psi\) & Obstacle function \\
\(f\) & Source term or forcing function \\
\(g\) & Boundary function \\
\(\mathcal{C}\) & Contact set where \(u=\psi\) \\
\(\Gamma(t)\) & Moving free boundary or interface \\
\(\Delta\) & Laplace operator \\
\(\mathcal{L}_{\mathrm{total}}\) & Total loss function \\
\(\mathcal{L}_{\mathrm{pde}}\) & \gls{pde} residual loss \\
\(\mathcal{L}_{\mathrm{ic}}\) & Initial-condition loss \\
\(\mathcal{L}_{\mathrm{obs}}\) & Obstacle-constraint loss \\
\(\mathcal{L}_{\mathrm{comp}}\) & Complementarity loss \\
\bottomrule
\end{tabular}
\end{table}

\section{Hyperparameter Settings in the Numerical Experiments}
\label{appendix:settings-pde}
This appendix summarizes the problem settings, network architectures, collocation-point sampling, and optimization hyperparameters used in the numerical experiments. Unless otherwise stated, all models are trained with the Adam optimizer, and all residual derivatives are computed using automatic differentiation.
\subsection{Obstacle Problem Experiments}
We consider two stationary obstacle benchmarks: the \gls{eop} and the nonlinear \gls{plop}.
The numerical settings are listed in Table~\ref{tab:obstacle_problem_settings}.
\begin{table}[!htbp]
\centering
\caption{Problem and training settings for the obstacle-type \glspl{bvp}.}
\label{tab:obstacle_problem_settings}
\begin{tabular}{lcc}
\toprule
\textbf{Quantity} & \textbf{\gls{eop}} & \textbf{\gls{plop}} \\
\midrule
Domain & \([0,1]^2\) & \([0,2]^2\) \\
Operator & \(-\Delta u\) & \(-\nabla\cdot(|\nabla u|^{p-2}\nabla u)\) \\
Source/forcing & \(f=0\) & \(f=-1\) \\
Order \(p\) & -- & \(p=4\) \\
Obstacle & Eq.~\eqref{eq:appendix_eop_obstacle} & Eq.~\eqref{eq:appendix_plop_obstacle} \\
Exact solution & Eq.~\eqref{eq:appendix_eop_exact_solution} & Eq.~\eqref{eq:appendix_plop_exact_solution} \\
Boundary data & \(g=u_{\mathrm{exact}}\) & \(g=u_{\mathrm{exact}}\) \\
KAN widths & \((2,48,48,48,1)\) & \((2,48,48,48,1)\) \\
Spline grid size & \(12\) & \(15\) \\
Epochs & \(10{,}000\) & \(10{,}000\) \\
Learning rate & \(10^{-3}\) & \(10^{-3}\) \\
Interior points & \(2048\) & \(4096\) \\
Boundary points & \(1024\) & \(1024\) \\
\gls{pde} penalty weight \(\lambda_{\mathrm{pde}}\) & \(1\) & \(1\) \\
Obstacle penalty weight \(\lambda_{\mathrm{obs}}\) & \(20\) & \(25\) \\
Complementarity penalty weight \(\lambda_{\mathrm{comp}}\) & \(1\) & \(1\) \\
Scheduler step & \(2000\) & \(2000\) \\
Scheduler factor & \(0.2\) & \(0.5\) \\
Evaluation grid & \(256\times256\) & \(256\times256\) \\
\bottomrule
\end{tabular}
\end{table}
\FloatBarrier
\subsection{Stefan Problem Experiments}
\begin{table}[!htbp]
\centering
\caption{Problem and training settings for the two-dimensional one-phase \gls{sp}.}
\label{tab:stefan_problem_settings}
\begin{tabular}{lc}
\toprule
\textbf{Quantity} & \textbf{Setting} \\
\midrule
Space-time sampling box & \([0,2.25]\times[0,1]\times[0,1]\) \\
Governing equation & \(\partial_t u-\Delta u=0\) \\
Interface temperature condition & \(u(s(x_2,t),x_2,t)=0\) \\
Stefan free-boundary condition & Eq.~\eqref{eq:appendix_stefan_problem} \\
Exact temperature solution & Eq.~\eqref{eq:appendix_stefan_exact_solution} \\
Exact moving free boundary & Eq.~\eqref{eq:appendix_stefan_exact_interface} \\
Initial data & Eq.~\eqref{eq:appendix_stefan_initial_conditions} \\
Temperature-network widths & \((3,48,48,48,1)\) \\
Interface-network widths & \((2,32,32,1)\) \\
Spline grid size & \(15\) \\
Epochs & \(10{,}000\) \\
Learning rate & \(10^{-3}\) \\
Interior space-time points & \(2048\) \\
Initial-condition points & \(512\) \\
Free-boundary/interface points & \(512\) \\
Heat-equation weight \(\lambda_{\mathrm{heat}}\) & \(1\) \\
Initial-condition weight \(\lambda_{\mathrm{ic}}\) & \(1000\) \\
Interface-temperature weight \(\lambda_{\mathrm{int}}\) & \(100\) \\
Stefan-condition weight \(\lambda_{\mathrm{fb}}\) & \(100\) \\
Scheduler step & \(2000\) \\
Scheduler factor & \(0.5\) \\
Evaluation resolution & \(256 \times 256\) \\
\bottomrule
\end{tabular}
\end{table}
\pagebreak
\section{Comparative Study of Trainable Parameters and KAN Spline Basis Size}
\label{appendix:comparison_study_kan_parameters}

\begin{figure}[H]
    \centering

    \includegraphics[width=\linewidth]
    {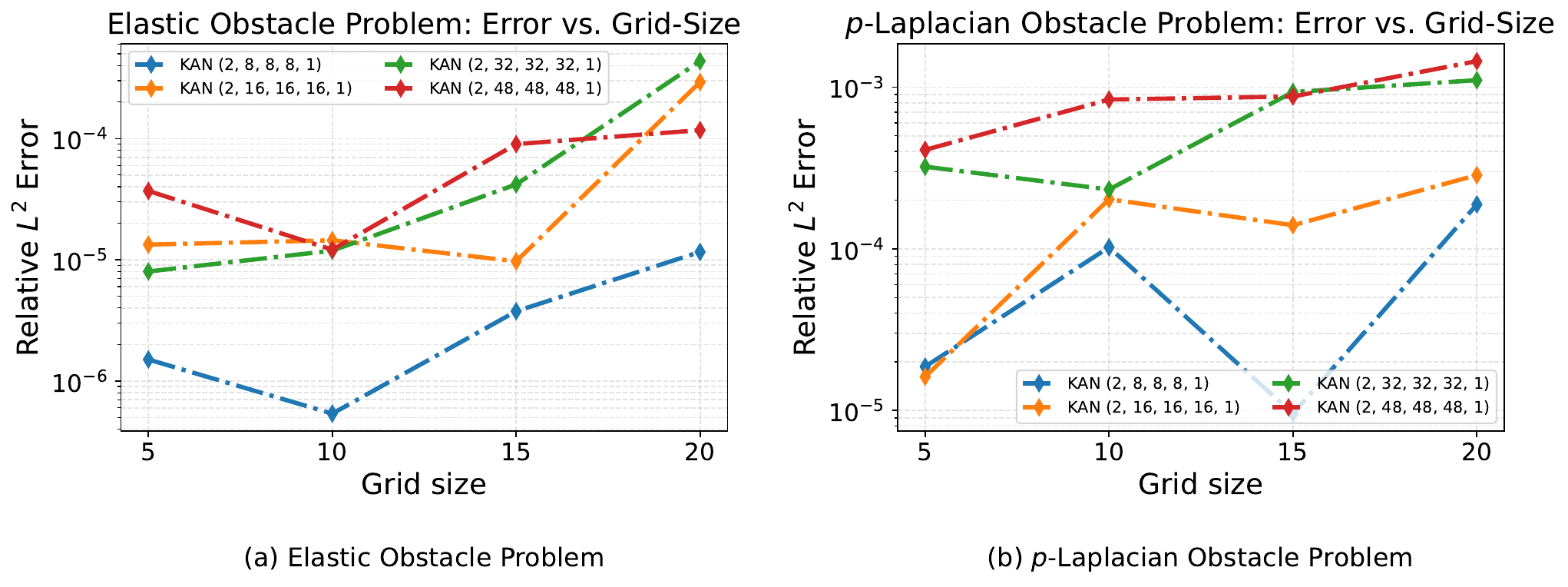}

    \vspace{0.5cm}

    \includegraphics[width=0.55\linewidth]
    {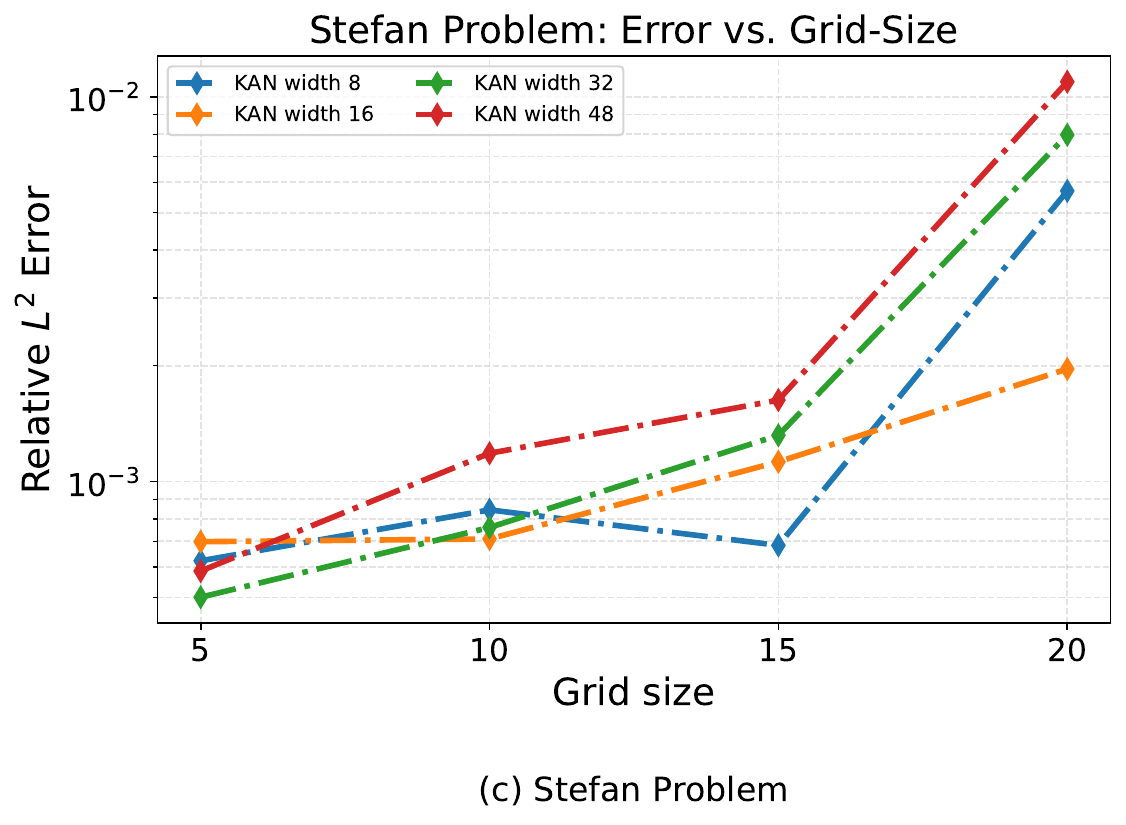}

    \caption{\textbf{Effect of \gls{kan} grid size and model width on obstacle-type and Stefan-type problems.}
    The top row reports the relative \(L^{2}\) error versus \gls{kan} grid size for four network architectures.
    Panel (a) corresponds to the elastic obstacle problem, while panel (b) corresponds to the nonlinear \(p\)-Laplacian obstacle problem.
    The bottom panel (c) reports the corresponding study for the Stefan problem.
    The results indicate that increasing the spline grid size does not necessarily lead to a monotonic reduction in error; rather, the approximation accuracy depends jointly on the spline grid resolution and network width.}

    \label{fig_params_gridsize_study}
\end{figure}

\end{document}